\documentclass[final]{article}

\usepackage{graphicx}
\usepackage{amsfonts}
\usepackage{amsmath}
\usepackage{amssymb}
\usepackage{fancyhdr}
\usepackage{titlesec}
\usepackage{indentfirst}
\usepackage{booktabs}
\usepackage{verbatim}
\usepackage{color}
\usepackage{latexsym}
\usepackage{amscd}
\usepackage{esvect}
\usepackage{graphicx}
\usepackage{upgreek}
\usepackage{subcaption}
\usepackage{caption}
\usepackage[page,header]{appendix}
\usepackage{titletoc}
\usepackage{mathrsfs}
\usepackage[numbers]{natbib}
\usepackage{float}
\usepackage{lipsum}
\usepackage[export]{adjustbox}
\usepackage[bookmarks=false]{hyperref}
\usepackage{epstopdf}

\numberwithin{equation}{section}
\newtheorem{theorem}{Theorem}[section]

\newtheorem{remark}[theorem]{Remark}

\allowdisplaybreaks

\def\bz{{\bf{z}}}

\def\bk{{\bf{k}}}
\def\bl{{\bf{l}}}

\begin{document}

\title{A stochastic asymptotic-preserving scheme for the bipolar semiconductor Boltzmann-Poisson system with random inputs 
and diffusive scalings 
\thanks{The author is supported by the funding DOE--Simulation Center for Runaway Electron Avoidance and Mitigation. }
}

\author{Liu Liu \footnote{The Institute for Computational Engineering and Sciences (ICES), University of Texas at Austin, 
Austin, TX 78712, USA (lliu@ices.utexas.edu).}}
\maketitle

\abstract{
In this paper, we study the bipolar Boltzmann-Poisson model, both for the deterministic system and the system with uncertainties, 
with asymptotic behavior leading to the drift diffusion-Poisson system as the Knudsen number goes to zero. 
The random inputs can arise from collision kernels, doping profile and initial data. 
We adopt a generalized polynomial chaos approach based stochastic Galerkin (gPC-SG) method. 
Sensitivity analysis is conducted using hypocoercivity theory for both the analytical solution and the gPC solution 
for a simpler model that ignores the electric field, and it gives their convergence toward the global Maxwellian exponentially in time. 
A formal proof of the stochastic asymptotic-preserving (s-AP) property and a {\it uniform} spectral convergence with error exponentially decaying in time 
in the random space of the scheme is given. Numerical experiments are conducted to validate the accuracy, 
efficiency and asymptotic properties of the proposed method. 
}


\section{Introduction}

Since kinetic equations are not first-principle physical equations, rather they often arise from mean field approximations of particle systems, hence there are inevitably modeling errors due to
incomplete knowledge of the interaction mechanism, imprecise measurement of the initial and boundary data, forcing terms, geometry, etc. These errors can  contribute uncertainties to the problems. Despite of intensive research at both theoretical and numerical levels, 
most researches are concerned with deterministic models and ignored uncertainties. Nevertheless, 
uncertainty quantification for kinetic equations, due to its importance in making reliable predications, calibrations and improvements of the kinetic models, deserves major attention from the research community.

To understand the propagation of the uncertainties and how they impact long-time behavior of the solution,
sensitivity and regularity analyses are crucial, since they allow us
to explore how sensitive the solution depends on random input parameters and to determine the convergence rate
of the numerical methods in the parameter space. In recent years one begins to see some activities in such studies, see for examples 
\cite{DesPer, JinLiuMa, Liu, QinWang, JinZhu, ShuJin, LiuJinUQ}.
At the numerical level, one of the popular UQ methods is the  generalized polynomial chaos method in the stochastic Galerkin (referred as gPC-SG)
framework \cite{GS, LMK, XiuBook}.
Compared with the classical Monte-Carlo method, the gPC-SG approach enjoys a spectral accuracy in the random space--provided the solution
is sufficiently regular in the space--while the Monte-Carlo method converges with only half-th order accuracy. 
{\color{blue}
As far as the non-intrusive stochastic collocation (SC) method is concerned, first the regularity analysis performed in this article is also useful for the accuracy analysis of SC methods. Second, there have been comparisons in terms of computational efficiencies between SG and SC for high dimensional problems; and there have been supporting cases the the SG methods are more efficient (see for example \cite{Xiu-Shen}).  
For the problem under study, it remains an interesting question to make such a comparison for high dimensional problems, but this is out of the scope of this article and could be an interesting future project. }
Recent studies of gPC-SG methods for kinetic equations and their behavior in various asymptotic regimes
are summarized in the review article \cite{HuReview}.

Kinetic equations play an important role in semiconductor device modeling \cite{MarkowichBook}. In such problems, the equations 
often have a diffusive scaling, characterized by the dimensionless Knudsen number $\varepsilon$, that leads asymptotically to the drift-diffusion equations as  
$\varepsilon$ goes to zero.
For multiscale problems in which $\varepsilon$ can vary in several orders of magnitude, the asymptotic-preserving (AP) schemes have proven to be effective and efficient to deal with different scales in a seamless way. An AP scheme switches between a micro solver and a macro one automatically, depending on the ratio of numerical parameters 
(mesh size, time step, etc.) over $\varepsilon$ \cite{Jin-AP-99, Jin-AP-Review, HJL-Review}. 
Just considering the transport of electrons in the conduction band, \cite{JinLorenzo} first introduced an AP scheme
for the semiconductor Boltzmann equation with an anisotropic collision operator, which is able to capture the correct diffusive behavior for the underresolved
numerical approximation. The scheme was further improved in \cite{Dengjia} with better stability condition. A higher-order scheme was constructed in \cite{Dimarco},
which improved the strict parabolic stability condition to a hyperbolic one. An efficient AP scheme in the high field regime was developed in
\cite{JinWang}. The authors in \cite{HuWang} further study the semiconductor Boltzmann equation with a two-scale stiff collision operators,
by taking into account different effects including the interactions between electrons and the lattice defects caused by ionized impurities \cite{Degond_ET}; 
they design and demonstrate the efficiency and accuracy of an asymptotic-preserving scheme that
leads to an energy-transport system as mean free path goes to zero at a discretized level.

For kinetic equations that contain random uncertainty, \cite{XZJ} first introduced the notion of
stochastic AP (s-AP), which was followed recently by many works successfully handling
the multiple scales for the kinetic equations with uncertainties \cite{Hu, JinLiu, MuLiu, JinLu}.
{\color{blue}s-AP scheme is introduced in the SG setting. It extends the idea from the deterministic AP 
methods to the stochastic case, which requires that as $\varepsilon\to 0$, the SG for the microscopic model with uncertainties 
automatically becomes a SG approximation for the limiting macroscopic stochastic equation. }

In this paper, we study the bipolar semiconductor system with random uncertainties, by taking into consideration
the generation-recombination process between electrons and holes \cite{Ansgar}.
The bipolar semiconductor Boltzmann-Poisson equations will be studied, and we
design and implement the gPC-SG scheme, with a formal proof of the s-AP property.
In order to analyze the convergence rate of the scheme, we use the hypocoercivity theory, which was well established in deterministic kinetic theory 
\cite{VillaniBook, DMS, CN, MB} and recently extended to study uncertain kinetic equations in the linear case \cite{QinWang} and nonlinear ones \cite{JinZhu, ShuJin, LiuJinUQ}.  By ignoring the self-consistent electric potential and using the hypocoercivity analysis done in \cite{MB, LiuJinUQ}, 
we obtain an exponential decay in time of the random solutions to the (deterministic) global equilibrium,
and uniform spectral convergence with an exponential decay in time of the numerical error of the gPC-SG method for the underlying system with uncertainties,
under suitable assumptions on the gPC polynomials and the random inputs.  To our knowledge,
this is the first study of AP and s-AP schemes for bipolar semiconductor Boltzmann system, in both deterministic and uncertain cases.

This paper is organized as the following. Section \ref{sec:2} gives an introduction of 
the bipolar Boltzmann-Poisson model, followed by a derivation of the limiting drift-diffusion equations. 
Section \ref{sec:3} discusses the AP scheme for the deterministic problem. A s-AP scheme in the gPC-SG framework 
for the bipolar model with random inputs will be studied and verified in section \ref{sec:4}. A convergence rate analysis for both the analytical solution 
and the gPC solution for a simpler model (without electric field) will also be conducted in section \ref{sec:4}. 
In section \ref{sec:6}, we present several numerical examples for both the deterministic problem and the model with uncertainties, to illustrate the efficiency, accuracy and s-AP properties of the proposed 
scheme. Finally, the paper is concluded in section \ref{sec:7}.  

\section{The bipolar semiconductor Boltzmann-Poisson system}
\label{sec:2}
In semiconductor devices, electrical currents originate from the transport of electrons and holes. 
$f_n(x,v,t)$, $f_p(x,v,t)$ represent the existence probability of an electron and a hole, respectively, at 
position $x\in\mathbb R^d$, with the velocity $v\in\mathbb R^d$, where $d$ is the dimension, at time $t\geq 0$. 
The Boltzmann equations that give the evolution of the distribution functions for them are written by (\cite{Ansgar, Poupaud})
\begin{eqnarray}
&\label{model_a}\displaystyle\epsilon \partial_t f_n + (v\cdot\nabla_x f_n -E\cdot \nabla_v f_n) =\frac{1}{\epsilon}Q_n(f_n)+\epsilon I_n(f_n, f_p), \\[6pt]
&\label{model_b}\displaystyle\epsilon \partial_t f_p +(\beta v\cdot\nabla_x f_p + E\cdot \nabla_v f_p)=\frac{1}{\epsilon}Q_p(f_p)+\epsilon I_p(f_n, f_p),  \\[6pt]
&\label{model_c} \displaystyle\gamma\, \Delta_x \Phi = n-p-C(x),  \qquad  E=-\nabla_x \Phi. 
\end{eqnarray}
where $\beta=m_{e}^{\ast}/m_{h}^{\ast}$ is the ratio of the effective masses of electrons and holes, which we consider it a constant. 
$\Phi=\Phi(t,x)$ represents the electric potential, $E=E(t,x)$ is the self-consistent electric field given by the Poisson equation (\ref{model_c}). 
$\gamma$ is some scaled Debye length, $C(x)$ is the doping profile. 
The densities of the electron and the hole is given by 
$$n=\int_{\mathbb R^d} f_n\, dv, \qquad p=\int_{\mathbb R^d} f_p\, dv. $$ 

Under the low density approximation, the linear collision operators are given by 
\begin{equation} Q_i (f_i)=\int_{\mathbb R^d} \sigma_i(x, v,w)(M_i(v)f_i(w)-M_i(w)f_i(v))dw,  \qquad i=n, \, p\,,  \label{Qi}\end{equation}
with 
\begin{equation} M_n(v)=\frac{1}{(2\pi)^{d/2}}e^{-|v|^2/2}\,,  \qquad M_p(v)=\frac{1}{(2\pi/\beta)^{d/2}}e^{-\beta|v|^2/2}\,. \label{Max}\end{equation}
being the normalized Maxwellian distribution of the electrons and holes. 
The anisotropic scattering kernel $\sigma_i$ for electrons and holes respectively are rotationally invariant and satisfies
\begin{equation}\label{sigma}\sigma_i(x,v,w)=\sigma_i(x,w,v)>0, \qquad i =n, \, p\,. \end{equation}

The process of {\it generation} of an electron-hole pair is that an electron moves from the valence band to the conduction band, 
leaving a hole behind it in the valence band. The inverse process of an electron moving from the conduction to the valence band is termed
the {\it recombination} of an electron-hole pair. See the following figure for an explanation \cite{Ansgar}: 
\begin{figure}[H]
\includegraphics[width=0.55\textwidth]{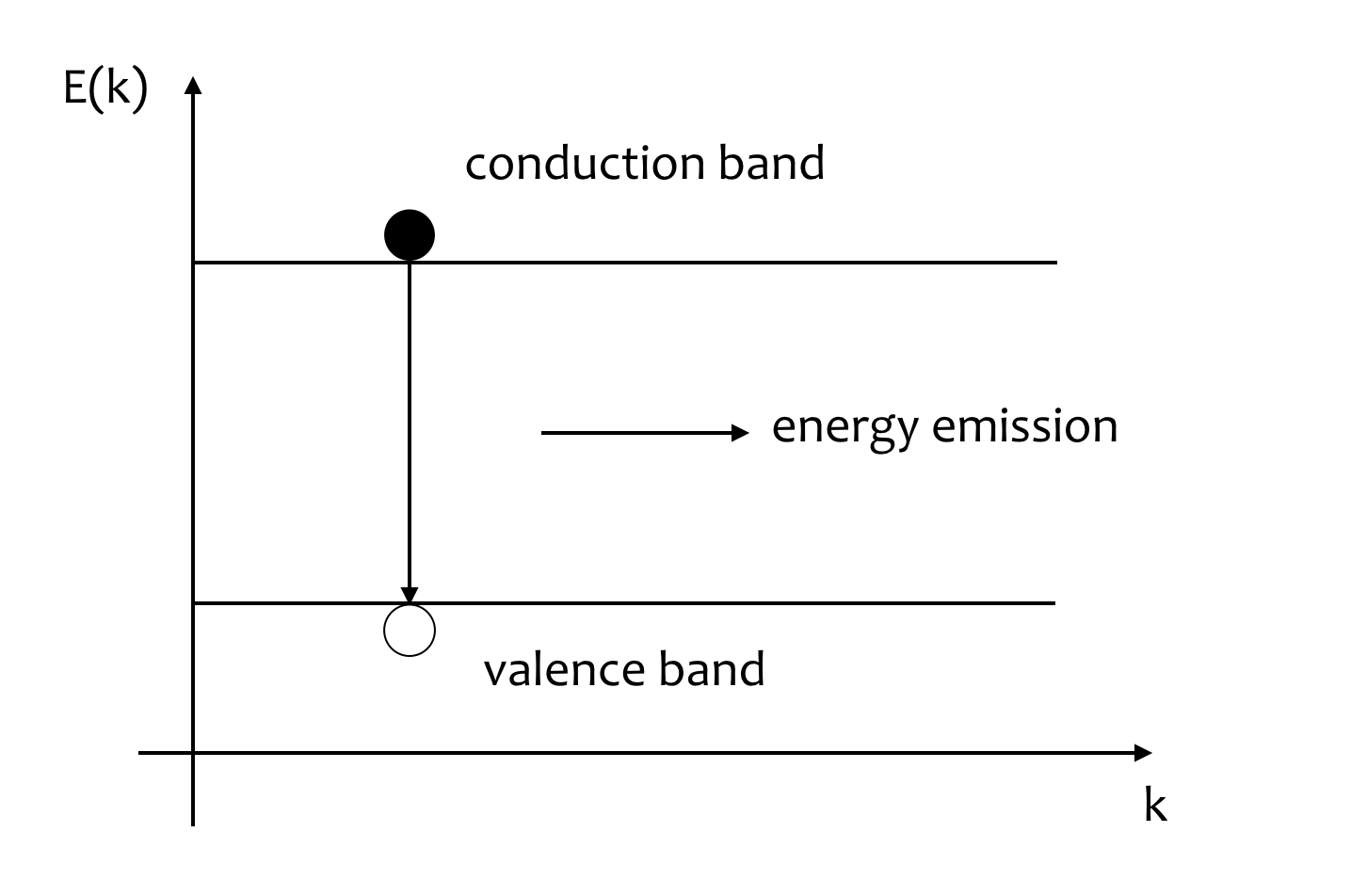}
\includegraphics[width=0.53\textwidth]{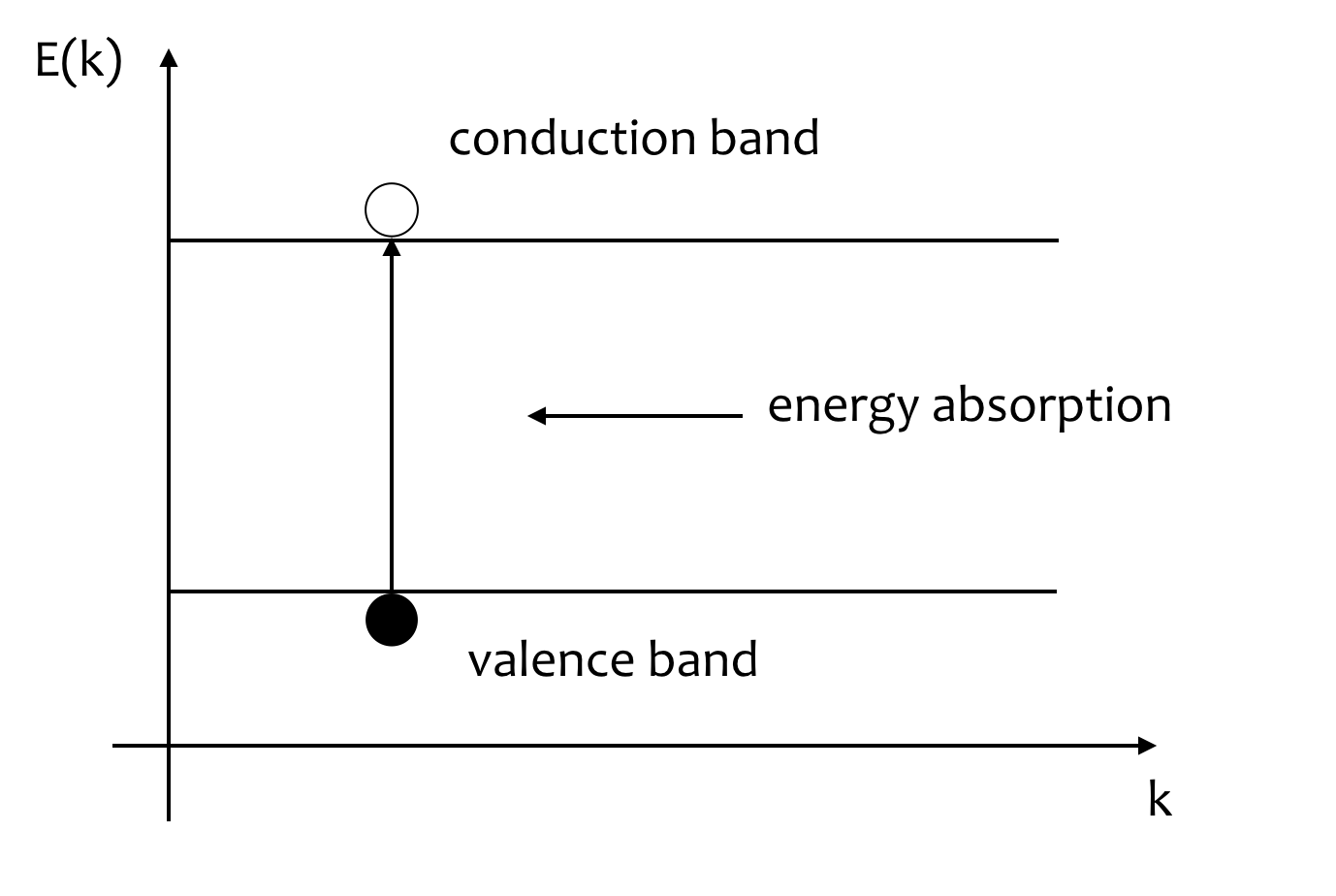}
\caption{A recombination-generation process}
\end{figure}

The recombination-generation operators are given by \cite{Ansgar, Poupaud}
\begin{eqnarray}
&\label{RG0}\displaystyle I_n(f_n, f_p)=\int_{\mathbb R^d} \sigma_I(x,v,w)\left[M_n(v)-M_p(w)f_n(v)f_p(w)\right]dw, \\[4pt]
&\label{RG}\displaystyle I_p(f_n, f_p)=\int_{\mathbb R^d} \sigma_I(x,w,v)\left[M_n(w)-M_p(v)f_n(w)f_p(v)\right]dw, 
\end{eqnarray}
where $\sigma_I$ is the generation-recombination kernel and is also rotationally invariant, as given in (\ref{sigma}). 
The collision frequency for electrons and holes is given by 
\begin{equation}\label{freq1} \lambda_i(x,v)=\int_{\mathbb R^d}\sigma_i(x,v,w)M_i(w)dw,  \qquad i=n, \, p\,. 
\end{equation}

The author in \cite{Poupaud} has proved the existence and uniqueness of smooth solutions of the system (\ref{model_a})--(\ref{model_c}). 

{\color{red} {\bf Remark.}
We give some explanations for the derivation of $I_n$, $I_p$ that model the generation and recombination processes:
\begin{equation}
\label{R:I_n} I_n(f_n, f_p)=\int_{\mathbb R^d} \sigma_I(x,v,w)\left[M_n(v)(1-f_n(v))(1-f_p(w)) - M_p(w) f_n(v) f_p(w)\right]dw. 
\end{equation}
The first term in the integral $I_n$ in (\ref{R:I_n}) represents the probability of creation of an electron at the coordinates $(x,v)$ 
and a hole at $(x,w)$; the second term in the integral represents the probability of recombination of an electron-hole pair. 
Due to the hypothesis of low density, i.e., $f_n, \, f_p \ll 1$, 
the terms $(1-f_n(v))$, $(1-f_p(w))$ tend to be $1$, then one gets $I_n$ defined in (\ref{RG0}). Similarly for $I_p$ given in (\ref{RG}).  

The recombination-generation effects are not negligible and crucial in many physics applications such as bipolar transistors, solar cells, LEDs and semiconductor lasers. Take solar cells as an example. Their mechanism is composed of several steps, that is, 
electron-hole pair generation by absorption of light in semiconductors, 
separation of electron-hole pairs by built-in potential, electron-hole recombination, etc. 
Understanding the recombination-generation processes is important and could help us improve the efficiency of solar cells \cite{solar}. 
}

\subsection{Bipolar drift-diffusion equations}
The relaxation time for collision and generation-recombination process has the relation: $\tau_{\text{col}}\gg \tau_{\text{gen}}$. 
Indeed, the typical time scale of collisions is $10^{-12}$s whereas the typical time for recombination-generation effects is $10^{-9}$s. 
In this section, the drift-diffusion equations are derived under the assumption that collisions occur on a much faster time scale than
recombination-generation processes. 

First, let us recall the following properties for the collision operators $Q_i$, for $i=n, p$\,, as discussed in \cite{Ansgar, Poupaud}. \\ 
\noindent (i) The kernel of $Q_i$ is spanned by $M_i$. 

\noindent(ii) $Q_i(f)=g$ has a solution if and only if $\int_{\mathbb R^d} g\, dv=0$. 

\noindent(iii) The equations $$Q_n(h_n)=vM_n(v),  \qquad Q_p(h_p)=\beta vM_p(v)$$
 have solutions $h_n$, $h_p$ with the property that there exist $\mu_{0,n}$, $\mu_{0,p}\geq 0$ satisfying
\begin{equation}\label{mu}\int_{\mathbb R^d} v\otimes h_n\, dv =-\mu_{0,n}I, \qquad 
\int_{\mathbb R^d} \beta v\otimes h_p\, dv = -\mu_{0,p}I, \end{equation}
where $I \in\mathbb R^{d\times d}$ is the identity matrix. 

Let $(f_n^{\epsilon}, f_p^{\epsilon}, E^{\epsilon})$ be a solution of (\ref{model_a})--(\ref{model_c}). 
As $\epsilon\to 0$ in (\ref{model_a}) and (\ref{model_b}), then
$$Q_n(f_n)=0, \qquad Q_p(f_p)=0, $$
where $\displaystyle f_i = \lim_{\epsilon\to 0} f_i^{\epsilon}$. Thus $f_n=nM_n$ and $f_p=pM_p$
by property (i). 

Inserting the Chapman-Enskog expansions 
\begin{equation}\label{Enskog} f_n^{\epsilon}=nM_n+\epsilon g_n^{\epsilon},
 \qquad f_p^{\epsilon}=pM_p+\epsilon g_p^{\epsilon} 
\end{equation}
into the Boltzmann equations (\ref{model_a}), one has
\begin{eqnarray*}
&\displaystyle \epsilon (n {\color{red}M_{n}}+\epsilon g_n^{\epsilon}) + (v\cdot\nabla_x(n{\color{red} M_{n}})-E\cdot\nabla_v(n 
{\color{red}M_{n}})) \\[6pt]
&\displaystyle \qquad\qquad\qquad\qquad\qquad\qquad\qquad + \epsilon (v\cdot\nabla_x g_n^{\epsilon}-E\cdot\nabla_v g_n^{\epsilon}) = Q_n(g_n^{\epsilon}) + \epsilon I_n(f_n^{\epsilon}, f_p^{\epsilon}). 
\end{eqnarray*}
The limit $\epsilon\to 0$ yields  
\begin{equation}\displaystyle Q_n(g_n)=(\nabla_x n+nE)\cdot v {\color{red} M_{n}}. \label{Qn_limit} \end{equation}
Similarly, inserting the expansion (\ref{Enskog}) into (\ref{model_b}), one gets
\begin{equation}\displaystyle Q_p(g_p)=\beta (\nabla_x p-pE)\cdot v {\color{red}M_{p}}. \label{Qp_limit}\end{equation}
where $\displaystyle g_i=\lim_{\epsilon\to 0}g_i^{\epsilon}$, $i=n,p$. 
By property (iii), solutions of (\ref{Qn_limit}) and (\ref{Qp_limit}) are
$$g_n=\frac{J_n}{\mu_{0,n}}\cdot h_n +c_n M_n, \qquad g_p=-\frac{J_p}{\mu_{0,p}}\cdot h_p +c_p M_p, $$
for some constants $c_n$, $c_p$, with $J_n$, $J_p$ defined by
$$ J_n=\mu_{0,n}(\nabla_x n+nE), \qquad J_p=-\mu_{0,p}(\nabla_x p-pE). $$ 
Thus 
\begin{equation}\langle v g_n\rangle =-J_n, \qquad \beta\langle v g_p \rangle =J_p, \label{vg}\end{equation}
where $\displaystyle\langle\, \cdot\, \rangle = \int_{\mathbb R^d} dv. $

Insert the Chapman-Enskog expansions (\ref{Enskog}) into  (\ref{model_a})--(\ref{model_b}) and integrate the velocity on both sides, then 
\begin{equation}
\begin{aligned}
\partial_t \langle nM\rangle +\epsilon\partial_t \langle g_n^{\epsilon}\rangle +\nabla_x\cdot\langle vg_n^{\epsilon}\rangle 
&=\langle I_n(nM_n+\epsilon\, g_n^{\epsilon}, \, pM_p+\epsilon\, g_p^{\epsilon})\rangle, \\[6pt]
\partial_t \langle pM\rangle +\epsilon\partial_t \langle g_p^{\epsilon}\rangle + \beta\nabla_x\cdot\langle vg_p^{\epsilon}\rangle
&=\langle I_p(nM_n+\epsilon\, g_n^{\epsilon}, \, pM_p+\epsilon\, g_p^{\epsilon})\rangle. 
\end{aligned}
\end{equation}
As $\epsilon\to 0$, by (\ref{vg}), one has 
\begin{eqnarray}
&\label{limit_a}\displaystyle \partial_t n -\nabla_x\cdot J_n =\langle I_n(nM_n, \, pM_p) \rangle, \\[6pt]
&\label{limit_b}\displaystyle  \partial_t p +\nabla_x \cdot J_p=\langle I_p(nM_n,\,  pM_p)\rangle. 
\end{eqnarray}
Denote $R(n,p)=\langle I_n(nM_n, pM_p) \rangle$, then 
\begin{align}
&\displaystyle R(n,p)= \int_{\mathbb R^d}\int_{\mathbb R^d}\sigma_I(x,v,w)M_n(v)dwdv - np \int_{\mathbb R^d}\int_{\mathbb R^d}
\sigma_I(x,v,w)M_p^2(w)M_n(v)dwdv\notag\\[4pt]
&\label{R_def}\displaystyle\qquad\quad := A(x) - np B(x), 
\end{align}
where we define 
$$A(x)=\int_{\mathbb R^d}\int_{\mathbb R^d}\sigma_I(x,v,w)M_n(v)dwdv, \qquad B(x)=\int_{\mathbb R^d}\int_{\mathbb R^d}
\sigma_I(x,v,w)M_p^2(w)M_n(v)dwdv. $$
Also note that $\langle I_n(nM_n, pM_p) \rangle=\langle I_p(nM_n, pM_p) \rangle$. 
The bipolar drift-diffusion Poisson system is given below: 
\\[2pt]

{\bf{\large Bipolar drift-diffusion equations}}
\begin{equation}
\begin{aligned}
\partial_t n -\nabla_x\cdot J_n &= R(n,p),  \qquad\qquad J_n=\mu_{0,n}\, (\nabla_x n+nE), \\[6pt]
\partial_t p -\nabla_x \cdot J_p &= R(n,p), \qquad\qquad J_p=\mu_{0,p}\, (\nabla_x p-pE),  \\[6pt]
-\gamma\, \nabla_x E &= n-p-C(x), \qquad x\in \mathbb R^d, 
\end{aligned}
\label{diff1}
\end{equation}
with $R(n,p)$ defined in (\ref{R_def}). 
\\[2pt]

{\color{red} {\bf Remark.} We list below some major differences and 
numerical difficulties compared with the single-species semiconductor Boltzmann equation studied in \cite{JinLorenzo}. 
We first recall the model equation (2.1) in \cite{JinLorenzo}, 
\begin{equation}\label{model-ref} \varepsilon \partial_t f + {\bf v}\cdot \nabla_{{\bf x}}f - \frac{q}{m}{\bf E}\cdot\nabla_{{\bf v}}f = \frac{1}{\varepsilon}{\bf Q}(f) 
+ \varepsilon G. \end{equation}
There $G=G(t, {\bf x}, {\bf v})$ is a source term that models the generation-recombination process. One can see that 
$G$ is not a function of $f$, thus the model studied in \cite{JinLorenzo} is linear, and only constant functions $G$ are considered in their numerical tests. 

In our model systems under study, $I_n(f_n, f_p)$ and $I_p(f_n, f_p)$ on the right-hand-side of (\ref{model_a})--(\ref{model_b})
model the generation and recombination of an electron-hole pair. Defined in (\ref{RG}), $I_n$, $I_p$ are non-linear integral operators in $f_n$, $f_p$ and 
are much more complicated than $G=G(t, {\bf x}, {\bf v})$ considered in \cite{JinLorenzo}! 
In fact, equations (\ref{model_a})--(\ref{model_b}) that describe the evolution of the distribution functions for electrons and holes are {\it coupled} through these non-linear integral operators, which is accounted for the {\it major difference} compared to the single-species model. 
As $\varepsilon\to 0$, the limiting system--bipolar drift-diffusion equations given in (\ref{diff1}) are also different from the drift-diffusion equation for the 
single-species, with the non-linear term $R(n,p)$ on the right-hand-side. 
Even for the deterministic bipolar model, it is not a trivial extension of the numerical method 
developed in \cite{JinLorenzo}. We would like to emphasize that 
this project is the {\it first} study of AP and s-AP schemes for bipolar semiconductor Boltzmann--Poisson system, in both deterministic and uncertain settings. 
}

\section{Parity equations and diffusive relaxation system}
\label{sec:3}
\subsection{Even- and Odd- Parity Equations}
Consider the one-dimensional velocity space $v\in\mathbb R$. Denote $f_1=f_n$, $f_2=f_p$, $\rho_1=n$, $\rho_2=p$ 
and rewrite the system (\ref{model_a}) as
\begin{align}
\label{model_1}\partial_t f_1 + \frac{1}{\epsilon}(v\cdot\nabla_x f_1-E\cdot\nabla_v f_1) &=\frac{1}{\epsilon^2}Q_1(f_1) + I_1(f_1,f_2),  \\[6pt]
\label{model_2}\partial_t f_2 + \frac{1}{\epsilon}(\beta v\cdot\nabla_x f_2+E\cdot\nabla_v f_2)&=\frac{1}{\epsilon^2}Q_2(f_2) + I_2(f_1,f_2), \\[6pt]
\label{model_3}\gamma\, \nabla_x E & =\rho_1-\rho_2-C(x). 
\end{align}

We will use the even- and odd- parities formulation, which is an effective vehicle to derive asymptotic-preserving scheme 
for linear transport equation \cite{JPT2} and one-component semiconductor Boltzmann equation \cite{JinLorenzo}. 
First, introduce the even parities $r_i$ and the odd parities $j_i$, for $i=1,2$, 
\begin{equation}
\begin{aligned}
r_i(t,x,v) & =\frac{1}{2}\left[f_i(t,x,v)+f_i(t,x,-v)\right], \\[6pt]
j_i(t,x,v)  &=\frac{1}{2\epsilon}\left[f_i(t,x,v)-f_i(t,x,-v)\right]. 
\end{aligned}
\label{RJ}
\end{equation}
Split (\ref{model_1}) and (\ref{model_2}) respectively into two equations, one for $v>0$ and one for $-v$, then 
\begin{equation}
\begin{aligned}
\partial_t f_i + \frac{1}{\epsilon}(s_i v\cdot\nabla_x f_i \mp  E\cdot\nabla_v f_i)& =\frac{1}{\epsilon^2}Q_i(f_i)(v) + I_i(f_1,f_2)(v), \\[6pt]
\partial_t f_i - \frac{1}{\epsilon}(s_i v\cdot\nabla_x f_i \pm E\cdot\nabla_v f_i) &=\frac{1}{\epsilon^2}Q_i(f_i)(-v) + I_i(f_1,f_2)(-v). 
\end{aligned}
\label{EO}
\end{equation}
where $s_1=1$, $s_2=\beta$. 
(A notation remark: in the first equation, $i=1$ corresponds to $- E\cdot\nabla_v f_1$ and 
$i=2$ corresponds to $E\cdot\nabla_v f_2$\, ; in the second equation, $i=1$ corresponds to $E\cdot\nabla_v f_1$ and 
$i=2$ corresponds to $-E\cdot\nabla_v f_2$\,.)

Adding (and multiplying by $1/2$), subtracting (and multiplying by $1/2\epsilon$) the two equations in (\ref{EO}), for $i=1$, $2$, respectively,
one gets
\begin{equation}
\begin{aligned}
\partial_t r_1 + v\cdot\nabla_x j_1 - E\cdot\nabla_v j_1 &= \frac{1}{\epsilon^2}Q_1(r_1)+I_{\text{1,plus}}(r_1,r_2), \\[6pt]
\partial_t j_1 + \frac{1}{\epsilon^2}(v\cdot\nabla_x r_1 - E\cdot\nabla_v r_1) &=
-\frac{1}{\epsilon^2}\lambda_1\, j_1 + I_{\text{1,minus}}(r_2, j_1), 
\end{aligned}
\label{EO1}
\end{equation}
and
\begin{equation}
\begin{aligned}
\partial_t r_2 + v\cdot\nabla_x j_2 + E\cdot\nabla_v j_2 &= \frac{1}{\epsilon^2}Q_2(r_2)+I_{2,\text{plus}}(r_1,r_2),  \\[6pt]
\partial_t j_2 + \frac{1}{\epsilon^2}(\beta v\cdot\nabla_x r_2 + E\cdot\nabla_v r_2) &=
-\frac{1}{\epsilon^2}\lambda_2\, j_2 + I_{\text{2,minus}}(r_1,j_2), 
\end{aligned}
\label{EO2}
\end{equation}
where
\begin{equation}
\begin{aligned}
 I_{\text{1,plus}}(r_1,r_2) &=\frac{1}{2}\int_{\mathbb R}\left(\sigma_I(v,w) + \sigma_I(-v,w)\right) dw M_1(v) - \int_{\mathbb R}
 \sigma_I(v,w)r_2(w)M_2(w)dw\, r_1(v),   \\[6pt]
 I_{\text{2,plus}}(r_1,r_2) &= \frac{1}{2}\int_{\mathbb R}\left(\sigma_I(v,w) + \sigma_I(-v,w)\right)M_1(w)dw
-\int_{\mathbb R}\sigma_I(v,w)r_1(w)dw\, r_2(v)M_2(v),  \\[6pt]
 I_{\text{1,minus}}(r_2, j_1) &=\frac{1}{2}\int_{\mathbb R}\left(\sigma_I(v,w) -\sigma_I(-v,w)\right) dw M_1(v) - \epsilon\int_{\mathbb R} \sigma_I(v,w)r_2(w)M_2(w)dw \, j_1(v),  \\[6pt]
 I_{\text{2,minus}}(r_1,j_2) &=\frac{1}{2}\int_{\mathbb R}\left(\sigma_I(v,w) - \sigma_I(-v,w)\right)M_1(w)dw - \epsilon \int_{\mathbb R}\sigma_I(v,w)r_1(w)dw\,  j_2(v)M_2(v), 
\end{aligned}
\label{I_DEF}
\end{equation}
which is derived in Appendix (ii). 

The macroscopic variables $\rho_i$ and mean velocity $u_i$ can be expressed in terms of the new variables
$r_i$, $j_i$ ($i=1,2$), 
\begin{equation}
\begin{aligned}
\rho_i (t,x) &=\int_{\mathbb R} f_i(t,x,v) dv= \int_{\mathbb R}  r_i(t,x,v) dv, \\[4pt]
 u_i(t,x) &=\frac{1}{\epsilon\rho_i}\int_{\mathbb R} f_i(t,x,v)v\, dv=\frac{1}{\rho_i}\int j_i(t,x,v)v\, dv. 
\end{aligned}
\label{mac}
\end{equation}

\subsection{Diffusive relaxation system}
As was done in \cite{JinLiu, JinLorenzo, JPT2}, we rewrite the equations (\ref{EO1})--(\ref{EO2}) into the following diffusive relaxation system
\begin{equation}
\begin{aligned}
\partial_t r_i + v\cdot\nabla_x j_i \mp E\cdot\nabla_v j_i &= \frac{1}{\epsilon^2}Q_i(r_i)+I_{i,\text{plus}}\,,  \\[6pt]
\partial_t j_i + \phi (s_i v\cdot\nabla_x r_i \mp E\cdot\nabla_v r_i) 
&=-\frac{1}{\epsilon^2}\left[\lambda_i j_i +(1-\epsilon^2\phi)(s_i v\cdot\nabla_x r_i \mp E\cdot\nabla_v r_i)\right]+I_{i,\text{minus}}\,, 
\end{aligned}
\label{EO5}
\end{equation}
where $\phi=\phi(\epsilon)$ is a control parameter such that $0\leq\phi\leq 1/\epsilon^2$.
One simple choice of $\phi$ is 
$$\phi(\epsilon)=\min\left\{1,\, \frac{1}{\epsilon^2}\right\}. $$
A standard time splitting on the system (\ref{EO5}) consists of a relaxation step 
\begin{align}
&\displaystyle\label{relax_1}\partial_t r_i = \frac{1}{\epsilon^2}Q_i(r_i), \\[6pt]
&\displaystyle\label{relax_2}\partial_t j_i = -\frac{1}{\epsilon^2}\left[\lambda_i\, j_i +(1-\epsilon^2\phi)(s_i v\cdot\nabla_x r_i \mp E\cdot\nabla_v r_i)\right], 
\end{align}
and the transport step 
\begin{equation}
\begin{aligned}
\partial_t r_i + v\cdot\nabla_x j_i \mp E\cdot\nabla_v j_i &= I_{i,\text{plus}}\,,  \\[6pt]
\partial_t j_i + \phi\, (s_i v\cdot\nabla_x r_i \mp E\cdot\nabla_v r_i) &= I_{i,\text{minus}}\,. 
\end{aligned}
\label{trans_1}
\end{equation}
\\[2pt]
{\color{red} {\bf Remark.} 

We address the major numerical difficulties compared to the single-species problem studied in \cite{JinLorenzo}. 
With the non-linear integral operators $I_1$, $I_2$ in (\ref{EO}), in order to use the even-odd
decomposition method, extra effort is needed to deal with the non-linear terms. 
The linear transport terms on the left-hand-side and the linear collision terms on the right-hand-side remain linear after adding and subtracting 
of the two equations in (\ref{EO}). The difficulty is to derive what the non-linear integral operators become, 
namely, to write the non-linear terms (after the addition and subtraction operations) as functions with respect to one of the pairs in the set $\{r_1, r_2, j_1, j_2\}$. This calculation requires repeatedly use of change of variables, 
rotationally invariance property and symmetry of the collision kernel $\sigma_I$, and is shown clearly in the Appendix. 
Moreover, these non-linear operators $I_{i,\text{plus}}$, $I_{i, \text{minus}}$ increase 
the computational complexity of the gPC-SG method introduced in section \ref{subsec:gPC}, 
where we have tensor products of matrices and vectors there. 

We mention another major difference in numerical scheme compared with the one species semiconductor Boltzmann equation studied in \cite{JinLorenzo}. Indeed, for each species, the procedure of rewriting the equations (\ref{EO1})--(\ref{EO2}) into the diffusive relaxation system and adopting the first-order time-splitting is similar to \cite{JinLorenzo}, except that one needs to determine whether to put the non-linear terms 
$I_{i, \text{plus}}$, $I_{i, \text{minus}}$ ($i=1, 2$) on the right-hand-side of equations in the relaxation step or the transport step. We design the scheme to put them in the transport step so that the AP property is guaranteed. Furthermore, sAP property of the discretized gPC-SG scheme for the underlying system in the stochastic case is proved in section \ref{sec:sAP}. 
}

\subsection{A discretized asymptotic-preserving scheme}
\label{Det_scheme}
In the relaxation step (\ref{relax_1}), since the collision term is stiff, 
one needs to treat it implicitly. The generation-recombination term is non-stiff, so one can leave it explicitly. 
It is hard to invert the collision operator $Q_i$ generally (especially for the anisotropic case). 
In {\color{blue}\cite{JinLorenzo}}, a Wild sum based time relaxation scheme, first proposed in \cite{GPT}, was adopted to handle the stiffness in the collision term. 
In \cite{Dengjia}, a fully implicit scheme for one-component semiconductor Boltzmann equation in the diffusive regime 
in which the more convenient BGK penalization method of Filbet-Jin
\cite{Filbet-Jin} was developed. Here we also use this approach. 
We reformulate (\ref{relax_1}) into the following form
\begin{equation}\label{relax_1a}\partial_t r_i = \underbrace{\frac{1}{\epsilon^2}\left[Q_i(r_i)-P_i(r_i)\right]}_{\text{less stiff}}+
\underbrace{\frac{1}{\epsilon^2}P_i(r_i)}_{\text{stiff}}.  \end{equation}
The first term on the right hand side of (\ref{relax_1a}) is non-stiff, or less stiff and less dissipative compared to the second term, thus it can be discretized 
explicitly, which avoids inverting the operator $Q_i$. The second term on the right hand side of (\ref{relax_1a}) is stiff or dissipative, thus will be treated 
implicitly. 

The discretized scheme for the system (\ref{relax_1a}) and (\ref{relax_2}) is given by
\begin{align}
&\label{Dis_Relax1}\displaystyle\frac{r_i^{\ast}-r_i^n}{\Delta t} = \frac{1}{\epsilon^2} \left[Q_i(r_i^n)-P_i(r_i^n)\right]+ \frac{1}{\epsilon^2}P_i(r_i^{\ast}), \\[6pt]
&\label{Dis_Relax2}\displaystyle\frac{j_i^{\ast}-j_i^n}{\Delta t}=-\frac{1}{\epsilon^2}\left[\lambda_i\, j_i^{\ast} +(1-\epsilon^2\phi)(s_i v\cdot\nabla_x r_i^{\ast} \mp E^{\ast}\cdot\nabla_v r_i^{\ast})\right]. 
\end{align}
where $P_i$ is the BGK operator, which is a linear operator and is asymptotically close to the collision term $Q_i(f)$, and is given by 
\begin{equation}\label{BGK} P_i(r_i) = \eta_i (\rho_i M_i(v) - r_i), \end{equation}
where $\eta_i$ is some constant chosen as the maximum value of the Fr$\acute{e}$chet derivative $\nabla Q_i(r_i)$ \cite{Filbet-Jin}. 
In particular for the anisotropic semiconductor Boltzmann case, it is addressed in \cite{JinWang} that $\eta_i$ should be chosen 
to satisfy $\eta_i > \max_{v}\lambda_i(x,v)$ for $i=1, 2$, where $\lambda_i$ is the collision frequency defined in (\ref{freq1}). 
\\[6pt]

{\bf{\large{A Discretized Scheme:}}} \\

For notation simplicity, we describe the spatial discretization in one dimension. 
Consider the spatial domain $\Omega=[x_L, x_R]$ which is
partitioned into $N$ grid cells with a uniform mesh size $\Delta x=1/N$. 
Define the left boundary $x_L$ as $x_{1/2}$, right boundary $x_R$ as $x_{N+1/2}$, choose the spatial grid points $\displaystyle x_{i-1/2}=x_{1/2}+(i-1)\Delta x$, for $i=1, \cdots, N+1$. The $i$-th interior cell is $\displaystyle [x_{i-1/2},x_{i+1/2}]$, for $i=1,\cdots, N$, with the cell average at time level $t^n$ given by 
$$ U_i^n=\frac{1}{\Delta x}\int_{x_{i-1/2}}^{x_{i+1/2}}  U(t^n, x,v)\, dx.$$
The velocity discretization is performed using spectral approximation based on the Hermite polynomials, 
which is equivalent to the moment method. We refer the reader to \cite{Schmeiser, JinLorenzo} for details. 

The scheme can be implemented as follows. 
\begin{itemize}
\item Step 1. \quad
Update $\rho_i^{\ast}$ and $r_i^{\ast}$. \\ 
Integrate (\ref{Dis_Relax1}) over $v$, note that $\displaystyle\int_{\mathbb R} Q_i(r_i)\, dv=0$ and 
$\displaystyle\int_{\mathbb R} P_i(r_i)\, dv=0$, then 
\begin{equation}\rho_i^{\ast}=\rho^n. \label{rho_star}\end{equation}

Denote 
$$\theta_1^{(i)}=\frac{\Delta t}{\epsilon^2+\eta_i \Delta t}\,. $$
By (\ref{Dis_Relax1}), (\ref{BGK}) and (\ref{rho_star}), one can update $r_i^{\ast}$: 
\begin{equation} r_i^{\ast}=r_i^n + \theta_1^{(i)}\,Q_i(r_i^n). \end{equation} 

-- Step 1.1. 
One can use any Poisson solver such as the spectral method to solve for $\Phi$, 
$$-\gamma\, \Delta_x\Phi =\rho_1-\rho_2-C(x), $$
then update the electric field $E^{\ast}$ by using the equation $E=-\nabla_x\Phi$ and 
a second order spatial discretization. 

\item Step 2. \quad
Update $j_i^{\ast}$. \\
Denote 
$$\theta_2^{(i)}=\frac{\epsilon^2}{\epsilon^2+\lambda_i\Delta t}\,, \qquad
\theta_3^{(i)}=\frac{\Delta t\,(1-\epsilon^2\phi)}{\epsilon^2 + \lambda_i\Delta t}\,. $$
(\ref{Dis_Relax2}) can be solved explicitly since we already have $r_i^{\ast}$, 
\begin{equation}\label{step2} j_i^{\ast}=\theta_2^{(i)}\,j_i^n - \theta_3^{(i)}\,(s_i v\cdot\nabla_x r_i^{\ast} \mp E^{\ast}\cdot\nabla_v r_i^{\ast}), \end{equation}
where $s_1=1$ and $s_2=\beta$. The spatial derivative of $f$ that appears in (\ref{step2}) is approximated using central difference, which allows 
one to implement the scheme explicitly and guarantee a second-order accuracy. 

\item Step 3.  Update $r_i^{n+1}$, $j_i^{n+1}$ in the transport step. \\
For notation simplicity, we focus on the case $i=1$. 
To define the numerical fluxes we used the second-order upwind scheme (with slope limiter) in the spatial direction (\cite{JPT, Jin_Xin}). 
In the $x$-direction the Riemann invariants are 
$$U_1=\frac{1}{2}(r_1+\phi^{-\frac{1}{2}}j_1), \qquad\qquad  V_1=\frac{1}{2}(r_1-\phi^{-\frac{1}{2}}j_1), $$
which move with the characteristic speed $\pm\sqrt{\phi}$\,. The second-order upwind discretization of 
$r\pm\phi^{-\frac{1}{2}}j$ (drop the subscript $1$ in $r_1$, $j_1$) is given by 
\begin{align*}
&\displaystyle\frac{1}{2}(r+\phi^{-\frac{1}{2}}j)_{i+\frac{1}{2}}
=\frac{1}{2}(r+\phi^{-\frac{1}{2}}j)_{i} + \frac{\Delta x}{4}\mu_i^{+}, \\[6pt]
&\displaystyle\frac{1}{2}(r-\phi^{-\frac{1}{2}}j)_{i+\frac{1}{2}}
=\frac{1}{2}(r-\phi^{-\frac{1}{2}}j)_{i+1} - \frac{\Delta x}{4}\mu_{i+1}^{-}, 
\end{align*}
where $\mu_{i}^{\pm}$ are the slope limiters of $r\pm\phi^{-\frac{1}{2}}j$ on the $i$-th cell at $(\ast)$-th time step. 

For $v>0$, let $\tau=\sqrt{\phi}\, v\,\frac{\Delta t}{\Delta x}>0$, then 
\begin{align}
&\displaystyle r_i^{n+1}=(1-\tau)r_i^{\ast} +\frac{\tau}{2}(r_{i+1}^{\ast} + r_{i-1}^{\ast}) - \frac{\tau}{2\sqrt{\phi}}(j_{i+1}^{\ast} - j_{i-1}^{\ast})\notag \\[6pt]
&\label{r1}\displaystyle\qquad\quad+\frac{\tau}{4}\Delta x (-\mu_i^{+}-\mu_{i+1}^{-}+\mu_{i-1}^{+}+\mu_i^{-})\pm\Delta t E^{\ast}\cdot\nabla_v j_i^{\ast}
 + \Delta t\, I_{i, \text{plus}}^{\ast}\,,  \\[6pt]
&\displaystyle j_i^{n+1}=(1-\tau)j_i^{\ast} + \frac{\tau}{2}(j_{i+1}^{\ast}+j_{i-1}^{\ast})-\frac{\sqrt{\phi}\tau}{2}(r_{i+1}^{\ast} -r_{i-1}^{\ast}) \notag\\[6pt]
&\label{j1}\displaystyle\qquad\quad+\frac{\tau}{4}\sqrt{\phi}\Delta x (-\mu_i^{+}+\mu_{i+1}^{-}+\mu_{i-1}^{+}-\mu_i^{-})\pm\phi\, \Delta t E^{\ast}\cdot\nabla_v r_i^{\ast}
 + \Delta t\, I_{i, \text{minus}}^{\ast}\,. 
\end{align}
The slope limiter is defined by
$$\mu_i^{\pm}=\frac{1}{\Delta x}\left[\pm r_{i\pm1}+\phi^{-\frac{1}{2}}j_{i\pm 1}\mp r_i
-\phi^{-\frac{1}{2}}j_i\right]\psi(\theta_i^{\pm}), $$
with $$\theta_i^{\pm}=\left(\frac{r_i\pm\phi^{-\frac{1}{2}}j_i-r_{i-1}\mp\phi^{-\frac{1}{2}}j_{i-1}}{r_{i+1}\pm\phi^{-\frac{1}{2}}
j_{i+1}-r_i\mp\phi^{-\frac{1}{2}}j_i}\right)^{\pm}, $$
and $\psi$ is the particular slope limiter function. A simple minmod slope limiter is chosen here, 
$$\psi(\theta)=\max\{0,\min\{1,\theta\}\}. $$ 

To update $r_2^{n+1}$, $j_2^{n+1}$, one needs to change $\tau$ to $\tau=\sqrt{\phi\beta}\,v\, \frac{\Delta t}{\Delta x}$, 
and $\phi$ to $\phi\beta$ in (\ref{r1}), (\ref{j1}), except that the term
$\pm\phi\,\Delta t\, E^{\ast}\cdot\nabla_v r_i^{\ast}$ remains the same in (\ref{j1}). 
\end{itemize}

\begin{remark}
The velocity discretization is performed using the Hermite quadrature rule, see \cite{Klar, JinLorenzo, JinLiu}. 
We denote $N_v$ as the number of quadrature points used in the numerical tests. 
\end{remark}

\section{The model with random inputs}
\label{sec:4}
In this section, the two-band semiconductor system with random inputs is considered. The collision kernels 
describing the transition rate between the same-species collisions or the generation-recombination process between different species can be uncertain, 
due to incomplete knowledge of the interaction mechanism. The uncertainties may also come from inaccurate measurement of the initial data, 
boundary data, and the doping profile $C(x,\bz)$. (\ref{model_1})--(\ref{model_3}) with random inputs is given by
\begin{equation}
\begin{aligned}
\partial_t f_1 + \frac{1}{\epsilon}(v\cdot\nabla_x f_1-E\cdot\nabla_v f_1) &=\frac{1}{\epsilon^2}Q_1(f_1)(x,\bz)+ I_1(f_1,f_2)(x,\bz), \\[6pt] \partial_t f_2 + \frac{1}{\epsilon}(\beta v\cdot\nabla_x f_2+E\cdot\nabla_v f_2) &=\frac{1}{\epsilon^2}Q_2(f_2)(x,\bz) + I_2(f_1,f_2)(x,\bz), \\[6pt]
-\gamma \nabla_x E &=\rho_1-\rho_2-C(x,\bz), \\[6pt]
f_i(0,x,v,z)&=f_{i,\text{in}}(x,v,z). 
\end{aligned}
\label{UQ}
\end{equation}

\subsection{Regularity and local sensitivity results}
\label{Converg1}
Conducting the convergence rate analysis on system (\ref{UQ}) with a self-consistent potential is complicated and remains a future work. 
For a discussion of Vlasov-Poisson-Fokker-Planck system with random initial data and small scalings, see \cite{JinZhu}. 
In this section, we consider the following system without electric potential (and let the mass ratio $\beta=1$ for simplicity), 
\begin{equation}
\begin{aligned}
\partial_t f_i + \frac{1}{\epsilon}v\cdot\nabla_x f_i &=\frac{1}{\epsilon^2}Q_i(f_i)(x,z)+ I_i(f_1,f_2)(x,z), \\[6pt] 
f_i(0,x,v,z)&=f_{i,\text{in}}(x,v,z), \qquad i=1,2, \qquad z\in I_z \subset \mathbb R\,. 
\end{aligned}
\label{uq_1}
\end{equation}
We will use the hypocoercivity theory to prove the exponential convergence of the random solutions toward the (deterministic) 
global equilibrium, in addition to spectral accuracy and exponential decay of the numerical error of the gPC-SG method. 
{\color{red} This is an example that the framework studied in \cite{LiuJinUQ} for general class of collisional kinetic models with random inputs can be 
{\it generalized}. The main differences are: here we have a multi-species system; and the non-linear integral operators $I_1$, $I_2$ 
own a different scaling compared to that of the linear collision operators $Q_1$, $Q_2$.  }

Here is a brief review of the solution estimate in \cite{LiuJinUQ}:
\begin{align}
\label{BP}
\left\{
\begin{array}{l}
\displaystyle \partial_t f + \frac{1}{\varepsilon} v\cdot\nabla_x f = \frac{1}{\varepsilon^2}\mathcal C(f,f), \\[4pt]
\displaystyle  f(0,x,v,z) = f_{\text{in}}(x,v,z), 
\end{array}\right.
\end{align}
where we consider the incompressible Navier-Stokes or diffusion scaling. 
$\mathcal C$ is a general class of collision operators, both the collision kernels
and the initial data depend on the random variable $z\in I_z$, with $I_z$ a compact domain. 

Under the perturbative setting, $f$ should be a small perturbation of the global equilibrium (Maxwellian) $\mathcal M$:
\begin{equation}
 f=\mathcal M+ \varepsilon M h, \qquad \mathcal M=\frac{1}{(2\pi)^{\frac{d}{2}}}\, e^{-\frac{|v|^2}{2}}, 
\label{per-f}
\end{equation}
where $M=\sqrt{\mathcal M}$. Applying this $f$ into (\ref{BP}), then the fluctuation $h$ satisfies
\begin{equation}\partial_t h + \frac{1}{\varepsilon} v\cdot\nabla_x h
=\frac{1}{\varepsilon^2}\mathcal L(h)+\frac{1}{\varepsilon}\mathcal F(h,h), 
\label{INS-scaling}
\end{equation}
where $\mathcal L$ is the linearized (around $\mathcal M)$ collision operator,
and $\mathcal F$ is the nonlinear remainder.

\underline{\it{Notations: }}
For two multi-indices $j$ and $l$ in $\mathbb N^{d}$, define $$\partial_l^j = \partial/\partial v_j\, \partial/\partial x_l\,. $$
For derivatives in $z$, we use the notation $$\partial_z^{\alpha} h = \partial^{\alpha}h\,. $$
Denote $||\cdot||_{\Lambda}:= ||\, ||\cdot||_{\Lambda_v}\, ||_{L^2_x}$.
Define the Sobolev norms
\begin{eqnarray}
&&||h||_{H_{x,v}^s}^2 = \sum_{|j|+|l|\leq s}\, ||\partial_l^j h||_{L^2_{x,v}}^2\,, \qquad
 ||h||_{H_{x,v}^{s,r}}^2 = \sum_{|m|\leq r}\, ||\partial^m h||_{H_{x,v}^s}^2\,, \\
&&||h(x,v,\cdot)||_{H^{s}_{x,v}H_z^r}^2 = \int_{I_z}\, ||h||_{H_{x,v}^{s,r}}^2 \pi(z)dz, 
\end{eqnarray}
in addition to the $\sup$ norm in $z$ variable,
$$ ||h||_{H_{x,v}^{s,r} L_z^{\infty}}=\sup_{z\in I_z}\, ||h||_{H_{x,v}^{s,r}}\,. $$

The following estimates on $h$ and the spectral accuracy of the SG methods are proved in \cite{LiuJinUQ}:

\underline{\it{Result I: }}
Assume  $||h(0)||_{H_{x,v}^s L_z^{\infty}} \leq C_{I}$, if $h$ is a solution of (\ref{INS-scaling}) in
$H_{x,v}^s$ for all $z$, then
\begin{equation}\label{thm2_1}  ||h(t)||_{H_{x,v}^{s, r}L_z^{\infty}}\leq  C_{I}\, e^{-\tau_s t}\,, \qquad
||h(t)||_{H_{x,v}^{s} H_z^r} \leq  C_{I}\, e^{-\tau_s t}\,, 
 \end{equation}
where $C_I$, $\tau_s$ are positive constants independent of $\varepsilon$. 
 
It is shown in \cite{CN} {\color{red} that} the deterministic, linear relaxation model satisfies all the Assumptions H1--H4, by taking 
$||\cdot||_{\Lambda_v}=||\cdot||_{L^2_v}$, then $||\cdot||_{\Lambda}=||\cdot||_{L^2_{x,v}}$. 
Assumption H5 is also satisfied for the non-linear operator 
$I_1$, $I_2$, that is, for each $z\in I_z$, $\exists\, k_0\in\mathbb N$ and a constant $C>0$ such that $\forall\, k\geq k_0$, 
$$||I_{i}(h,h)||_{H^k_{x,v}}\leq C ||h||_{H^k_{x,v}}^2, $$
by using Sobolev embeddings and the Cauchy-Schwarz inequality, exactly the same as discussed in \cite{CN, MB}. 
The following assumptions on the collision kernels $\sigma_{i}$ ($i=1,2$) and $\sigma_{I}$ are needed:
\begin{equation} \label{assump1}
 |\partial_{z}^{k}\sigma_{i}(x,v,w,z)|\leq C_b, \qquad  |\partial_{z}^{k}\sigma_{I}(x,v,w,z)|\leq C_b^{\ast}, \qquad \forall\, k\geq 0\,. \end{equation}
Under these conditions, one can easily check that Assumptions H1--H5 given in \cite{LiuJinUQ} still hold 
when uncertainties are from collision kernels. 

Let $$f_i=\mathcal M + \varepsilon M h_i,  \qquad i=1, 2. $$
Plug it into the system (\ref{uq_1}), the perturbed solution $h_i$ satisfies
\begin{equation}\label{INS1}\partial_t h_i + \frac{1}{\varepsilon}v\cdot\nabla_x h_i =\frac{1}{\varepsilon^2}Q_i(h_i)   +  \varepsilon I_i(h_1, h_2). \end{equation}
Since the generation-recombination process has a weaker effect than the collision among particles, which leads to 
the non-linear operators $I_1$, $I_2$ owning a different scaling than the linear operators $Q_1$, $Q_2$. 
Note that whatever discussed in \cite{LiuJinUQ} for the scaled equation (\ref{INS-scaling}) remains valid for the problem we consider here, 
since the coefficient in front of the non-linear operator in (\ref{INS1}) and (\ref{INS-scaling}) has the relation: $\varepsilon < 1/\varepsilon$. 

Based on the proof of Lemma 3.1 in section 3 in \cite{LiuJinUQ}, as a corollary, it is obvious to check that the perturbed solution $h_i$ 
for the two-species system has the following estimate: 
\begin{equation}\label{hi}\frac{d}{dt} ||h_i||_{\mathcal H_{\epsilon_{\perp}}^{s,r}}^2 \leq \bigg[K_1 \sum_{i=1}^{2} ||h_i||_{H^{s,r}}^2 - K_2\bigg]
 \left(\sum_{i=1}^{2} ||h_i||_{H_{\Lambda}^{s,r}}^2 \right), \end{equation}
where the complicated definition of the norm $||\cdot||_{\mathcal H_{\epsilon_{\perp}}^s}$ is omitted, but one can check (2.20) in \cite{LiuJinUQ}. 
One just needs to know that 
$||\cdot||_{\mathcal H_{\epsilon_{\perp}}^s}$ is equivalent to $||\cdot||_{H^s}$, and that $||\cdot||_{H_{\Lambda}^{s,r}}=||\cdot||_{H^{s,r}}$ in our problem, 
then (\ref{hi}) becomes 
$$\frac{d}{dt}\left(\sum_{i=1}^{2} ||h_i||_{\mathcal H_{\epsilon_{\perp}}^{s,r}}^2\right)
\leq \bigg[K_3 \sum_{i=1}^{2}||h_i||_{\mathcal H_{\epsilon_{\perp}}^{s,r}}^2-K_2\bigg] \left(\sum_{i=1}^{2}||h_i||_{H^{s,r}}^2\right), $$
where $K_1, K_2, K_3$ are all constants independent of  $\varepsilon$ and $z$. 
 
If the initial data satisfies 
\begin{equation}
\label{IC} ||h_1(0)||_{\mathcal H_{\epsilon_{\perp}}^{s,r}}^2 + ||h_2(0)||_{\mathcal H_{\epsilon_{\perp}}^{s,r}}^2 
\leq \frac{K_2}{2 K_3}\,, 
\end{equation}
then $$\frac{d}{dt}\left(\sum_{i=1}^{2} ||h_i||_{\mathcal H_{\epsilon_{\perp}}^{s,r}}^2 \right) 
\leq {\color{red} -\frac{K_2}{2}}\left(\sum_{i=1}^{2} ||h_i||_{H^{s,r}}^2 \right) \leq 
- \widetilde C \left(\sum_{i=1}^{2} ||h_i||_{\mathcal H_{\epsilon_{\perp}}^{s,r}}^2\right). $$
The last inequality is because $H^{s}$ norm is equivalent to $\mathcal H_{\epsilon_{\perp}}^{s}$ norm. 

\begin{theorem}
\label{thm_f}
If the assumptions for the random kernels and the initial data--(\ref{assump1}) and (\ref{IC}) are satisfied, 
then solution of each species has the following estimate:
\begin{equation} ||h_i(t)||_{H_{x,v}^{s, r}L_z^{\infty}}\leq C_1\, e^{-\tau_1 t}\,, \qquad 
||h_i(t)||_{H_{x,v}^{s} H_z^{r}} \leq C_1\, e^{-\tau_1 t}\,, \,    i=1, 2, \end{equation}
where $C_1$, $\tau_1$ are constants independent of  $\varepsilon$ and $z$. 
\end{theorem}

This result shows that the random perturbation in both initial data and collision kernel will decay exponentially, 
and the random solutions $f_1(t,x,v,z)$, $f_2(t,x,v,z)$ will both converge exponentially in time to the deterministic global Maxwellian $\mathcal M$.
That is, the dependence on the random parameter $z$ of the two-band system is insensitive for long time. 

{\color{red} {\bf Remark.}
Thanks to the small $\mathcal O(\varepsilon)$ scaling of the non-linear integral terms $I_1$, $I_2$,  the analysis and conclusions presented in \cite{LiuJinUQ} can be extended here. Though compared to the previous work, 
where a complete framework for the kinetic equations with multiple scales and uncertainties and its gPC-SG systems has been well-established,  
the analysis conducted here is not as exquisite, yet it is a nice observation that the conclusions there can be adopted and generalized, 
since \cite{LiuJinUQ} does not mention directly the kinetic equation whose right-hand-side has a linear collision operator combined with a non-linear integral term and of different scalings. More importantly, this first attempt to study the bipolar semiconductor Boltzmann-Poisson system with random inputs from both numerical and analysis points of view may intrigue new directions of study.  
For example, conducting sensitivity analysis for the multi-species full Boltzmann equations with random inputs, which is more complicated and a non-trivial extension of the single-species problem studied in \cite{LiuJinUQ}. 
}

\subsection{A gPC-SG Method}
\label{subsec:gPC}
Let $\mathbb P_P^n$ be the space of the $n$-variate polynomials of degree less than or equal to $P\geq 1$, 
and recall that
$$ \text{dim}(\mathbb P_P^n)= \mbox{card}\{\bk \in \mathbb N ^n, |\bk|\leq P\}= \left(\begin{array}{c} n+P \\ P \end{array} \right):=K, $$
where we have denoted $\bk= (k_1,\dots, k_n)$ and $|\bk|=k_1+\dots+k_n$.
We consider the inner product
$$  \langle f, g\rangle_\pi =  \int_{I_{\bz}} f g\, \pi(\bz)d\bz, \quad \forall\, f, g \in L^2(\pi(\bz)d\bz),$$
where $L^2(\pi(\bz)d\bz)$ is the usual weighted Lebesgue space, and its associated norm is
$$ \|f\|_{L^2(\pi(\bz)d\bz)}^2 = \int_{I_{\bz}}f^2\, \pi(\bz)d\bz.$$
Consider a corresponding orthonormal basis $\{\psi_{\bk}(\bz)\}_{\bk\in \mathbb N ^n, \, |\bk|\leq P}$ of the space $\mathbb P_P^n$,  
where the degree of $\psi_\bk$ is $\mbox{deg}(\psi_\bk)= |\bk|$. In particular
\begin{equation*}
\langle \psi_\bk, \psi_\bl\rangle_\pi = \int_{I_{\bz}} \psi_\bk(\bz)\psi_\bl(\bz)\pi(\bz)d\bz=\delta_{\bk\bl}, \qquad |\bk|, \, |\bl|\leq P,
\end{equation*}
where $\delta_{\bk\bl}$ is the Kronecker symbol.
The commonly used pairs of $\{\psi_\bk(\bz)\}$ and $\pi(\bz)$ include Hermite-Gaussian, Legendre-uniform, Laguerre-Gamma, etc \cite{XiuBook, XiuKarn}. 
 
The SG method seeks the solution as a projection onto the space $\mathbb P_P^n$ (the set of $n$-variate 
orthonormal polynomials of degree up to $P\geq 1$), that is
\begin{equation}\label{soln u}
f(t,x,v,\bz) \approx f^K(t,x,v,\bz)=\sum_{k=1}^{K} \hat f_k(t,x,v)\psi_k(\bz). 
\end{equation}
From this approximation, one can easily compute statistical moments, such as
the mean and standard deviation, 
\begin{equation}\label{moments}
\mathbb{E}(f)\approx \hat f_1, \qquad
\text{SD}(f)\approx \big(\sum_{k=2}^{K}|\hat f_k|^2\big)^{1/2}\,. 
\end{equation}

By the gPC-SG approach, one inserts the ansatzes
\begin{equation}
\label{ans1}
\begin{aligned}
 f_i^K &=\sum_{k=1}^K \hat{(f_i)}_k \psi_k(\bz)={\bf \hat f}_{i}\cdot {\boldsymbol\psi}(\bz),  \qquad i=1, 2,  \\[4pt]
  E^K &=\sum_{k=1}^K \hat{E}_k \psi_k(\bz) = {\bf \hat E}\cdot {\boldsymbol\psi}(\bz)  
\end{aligned}
\end{equation}
into system (\ref{UQ}) and enforces the residual to be orthogonal to the polynomial space spanned by $\psi_\bk(\bz)$, then 
\begin{equation}
\begin{aligned}
\partial_t \hat{(f_1)}_k + \frac{1}{\epsilon}[v\cdot\nabla_x \hat{(f_1)}_k - \sum_i\sum_j \hat{E}_i \cdot \nabla_v \hat{(f_1)}_j\, G_{ijk}]
&=\frac{1}{\epsilon^2}({\bf Q}_1)_{k}({\bf \hat f}_1) + ({\bf I}_1)_{k}({\bf \hat f}_1, {\bf \hat f}_2),  \\[6pt]
\partial_t \hat{(f_2)}_k + \frac{1}{\epsilon}[\beta v\cdot\nabla_x \hat{(f_2)}_k + \sum_i\sum_j \hat{E}_i \cdot \nabla_v \hat{(f_2)}_j\, G_{ijk}]
&=\frac{1}{\epsilon^2}({\bf Q}_2)_{k}({\bf \hat f}_2) + ({\bf I}_2)_{k}({\bf \hat f}_1, {\bf \hat f}_2), \\[6pt]
- \gamma\nabla_x {\bf \hat E} &= (\hat{\rho}_1)_k - (\hat{\rho}_2)_k - L_k, 
\end{aligned}
\label{E_k}
\end{equation}
where 
\begin{align*}
&\displaystyle ({\bf Q}_i)_{k}({\bf \hat f}_{i}) = \int_{\mathbb R^d} (B_i(v,w))_k \left[M(v){\bf \hat f}_{i}(w) - M(w){\bf \hat f}_{i}(v)\right] dw, \qquad i=1,2,  \\[6pt]
&\displaystyle ({\bf I}_1)_{k}({\bf \hat f}_1, {\bf \hat f}_2) = \int_{\mathbb R^d} D_k(x,v,w)M_1(v)dw - \int_{\mathbb R^d} \sum_i \sum_j (\hat {f_1}(v))_i 
(\hat {f_2}(w))_j M_2(w) F_{ijk}(x,v,w) dw, \\[6pt]
&\displaystyle ({\bf I}_2)_{k}({\bf \hat f}_1, {\bf \hat f}_2) = \int_{\mathbb R^d} D_k(x,w,v)M_1(w)dw - \int_{\mathbb R^d} 
\sum_i \sum_j (\hat {f_1}(w))_i
(\hat {f_2}(v))_j M_2(v) F_{ijk}(x,v,w) dw, 
\end{align*}
with $(B_i)_k$ the $k$-th row of $K\times K$ matrix $(B_i)_{mn}$ ($i=1, 2$), given by 
\begin{equation}\label{B_matrix}
(B_i)_{mn}(x,v,w) = \int_{I_{\bz}}\sigma_i(x,v,w,\bz)\psi_m(\bz)\psi_n(\bz)\pi(\bz) d\bz. 
\end{equation}
The tensors $(G_{ijk})_{K\times K \times K}$, $(F_{ijk})_{K\times K \times K}$ and the vectors $(L_k)_{K\times 1}$, $(D_k)_{K\times 1}$ are defined by
\begin{equation}
\begin{aligned}
G_{ijk} &= \int_{I_{\bz}}\psi_i(\bz)\psi_j(\bz)\psi_k(\bz)\pi(\bz) d\bz, \\[6pt]
F_{ijk}(x,v,w) &=\int_{I_{\bz}}\sigma_I(x,v,w,\bz)\psi_i(\bz)\psi_j(\bz)\psi_k(\bz)\pi(\bz)d\bz, \\[6pt]
L_k(x) &=  \int_{I_{\bz}}C(x,\bz) \psi_k(\bz)\pi(\bz)d\bz, \\[6pt]
D_k(x,v,w) & = \int_{I_{\bz}} \sigma_I(x,v,w,\bz)\psi_k(\bz)\pi(\bz)d\bz. 
\end{aligned}
\end{equation}
\\[15pt]
{\bf\large{A convergence rate analysis}} \\

Here is a brief review of the gPC error estimate in \cite{LiuJinUQ} for the single species model: \\
\underline{\it{Result II: }}
Define the norm $$ ||h^e||_{H_{x,v}^{s} L_z^2}:= \int_{I_z}\, ||h^e||_{H_{x,v}^s}\, \pi(z)dz. $$
Under the technical conditions on the gPC polynomials:
\begin{equation}
\label{basis}
||\psi_k||_{L^{\infty}} \leq C k^p, \qquad \forall\, k,  
\end{equation}
we have
\begin{equation}
||h-h^K||_{H_{x,v}^{s}L_z^2} \leq C_{e}\,  \frac{e^{-\lambda t}}{K^r}\,, 
\end{equation}
with the constants $C_{e}, \,\lambda>0$ independent of $K$ and $\varepsilon$. 

We now give the main conclusion for the gPC error estimate for our two-band model: 
\begin{theorem}
\label{thm_gPC}
Under the assumption for the gPC polynomials (\ref{basis}) and the random kernels (\ref{assump1}), also assume that $\sigma_{I}$ is linear in $z$ with 
\begin{equation} 
\label{assump2} |\partial_{z}\sigma_{I}|\leq O(\varepsilon), 
\end{equation}
then  
\begin{equation} ||h_i-h_i^K||_{H_{x,v}^{s}L_z^2} \leq C_2\, \frac{e^{-\tau_2 t}}{K^r}\,, \qquad\text{for  } i=1, 2, \end{equation}
where $C_2$, $\tau_2$ are constants independent of  $\varepsilon$ and $z$. 
\end{theorem}
The proof of this theorem is really similar to \cite{LiuJinUQ} and we omit it here. Compared to \cite{LiuJinUQ}, 
one only needs to add up the estimates for $i=1$ and $i=2$, the same way as shown in the proof of Theorem \ref{thm_f}. 

To conclude, Theorem \ref{thm_gPC} gives a {\it uniform} spectral convergence of the SG method for the system (\ref{uq_1}), 
with convergence rate exponentially decaying in time, under suitable assumptions (\ref{assump1}), (\ref{basis}) and (\ref{assump2}). 
\\[20pt]
{\bf\large{The even-odd decomposition method}}\\

We use the even-odd decomposition and insert the ansatzes
$$ r_i^K = \sum_{k=1}^K \hat{(r_i)}_k \psi_k(\bz) ={\bf \hat r}_{i}\cdot {\boldsymbol\psi}(\bz), 
\qquad  j_i^K = \sum_{k=1}^K \hat{(j_i)}_k \psi_k(\bz) ={\bf\hat j}_{i}\cdot {\boldsymbol\psi}(\bz),  \qquad i=1, 2, $$
and $E^K$ in (\ref{ans1}) into systems (\ref{relax_1}) and (\ref{trans_1}). 
By the standard Galerkin projection, one gets the relaxation step 
\begin{align}
&\label{s_relax1}\displaystyle \partial_t (\hat {r_i})_k = \frac{1}{\epsilon^2} ({\bf Q}_i)_{k}({\bf \hat r}_{i}), \\[6pt]
&\label{s_relax2}\displaystyle \partial_t (\hat {j_i})_k = -\frac{1}{\epsilon^2} \left[(H_i)_k\, {\bf \hat j}_{i} + (1-\epsilon^2\Phi)(v\cdot\nabla_x (\hat r_i)_k \mp 
\sum_m\sum_n \hat E_m \cdot\nabla_v (\hat {r_i})_n G_{mnk})\right], 
\end{align} 
where $(H_i)_k$ is the $k$-th row of the matrix $(H_i)_{K\times K}$, given by 
$$(H_i)_{mn}(x,v) = \int_{I_{\bz}}\lambda_{i}(x,v,\bz)\psi_m(\bz)\psi_n(\bz)\pi(\bz)d\bz, $$
with the matrix $B_i$ given in (\ref{B_matrix}). The transport step is given by 
\begin{align}
\label{s_trans1}\displaystyle \partial_t (\hat r_i)_k + v\cdot\nabla_x (\hat j_i)_k  \mp \sum_m\sum_n \hat E_m \cdot \nabla_v (\hat j_i)_n\, G_{mnk} 
&= ({\bf I}_{\text{i,plus}})_k, \\[6pt]
\label{s_trans2}\displaystyle \partial_t (\hat j_i)_k + \Phi[s_i v\cdot \nabla_x (\hat r_i)_k \mp \sum_m \sum_n \hat E_m \cdot \nabla_v (\hat r_i)_n\, G_{mnk}]
&= ({\bf I}_{\text{i,minus}})_k, 
\end{align} 
where
\begin{equation*}
\begin{aligned}
({\bf I}_{\text{1,plus}})_k &= \frac{1}{2}M_1(v) J_k^a(x,v,w) - \int_{\mathbb R^d} \sum_m \sum_n (\hat r_1(v))_m (\hat r_2(w))_n\, 
M_2(w) F_{mnk}(x,v,w)dw, \\[6pt]
({\bf I}_{\text{2,plus}})_k  &= \frac{1}{2} J_k^c(x,v,w) -  \int_{\mathbb R^d} \sum_m \sum_n (\hat r_1(w))_m (\hat r_2(v))_n\, 
M_2(v) F_{mnk}(x,v,w)dw,  \\[6pt]
({\bf I}_{\text{1,minus}})_k &= \frac{1}{2}M_1(v) J_k^b(x,v,w) - \epsilon\int_{\mathbb R^d}  
\sum_m \sum_n (\hat j_1(v))_m (\hat r_2(w))_n\, 
M_2(w) F_{mnk}(x,v,w)dw, \\[6pt]
({\bf I}_{\text{2,minus}})_k &= \frac{1}{2} J_k^d(x,v,w) -\epsilon\int_{\mathbb R^d}\sum_m \sum_n (\hat r_1(w))_m (\hat j_2(v))_n\, 
M_2(w) F_{mnk}(x,v,w)dw, 
\end{aligned}
\end{equation*}
with $(J_k^a)_{K\times K}$, $(J_k^b)_{K\times K}$, $(J_k^c)_{K\times K}$ and $(J_k^d)_{K\times K}$ given by 
\begin{align*}
&\displaystyle J_k^a(x,v,w)= \int_{I_{\bz}}\int_{\mathbb R^d} \big(\sigma_I(x,v,w,\bz)+\sigma_I(x,-v,w,\bz)\big) dw\, \psi_k(\bz)\pi(\bz)d\bz, \\[6pt]
&\displaystyle J_k^b(x,v,w)= \int_{I_{\bz}}\int_{\mathbb R^d} \big(\sigma_I(x,v,w,\bz)-\sigma_I(x,-v,w,\bz)\big) dw\, \psi_k(\bz)\pi(\bz)d\bz, \\[6pt]
&\displaystyle J_k^c(x,v,w)= \int_{I_{\bz}}\int_{\mathbb R^d} \big(\sigma_I(x,v,w,\bz)+\sigma_I(x,-v,w,\bz)\big) M_1(w) dw\, \psi_k(\bz)\pi(\bz)d\bz, \\[6pt]
&\displaystyle J_k^d(x,v,w)= \int_{I_{\bz}}\int_{\mathbb R^d} \big(\sigma_I(x,v,w,\bz)-\sigma_I(x,-v,w,\bz)\big) M_1(w) dw\, \psi_k(\bz)\pi(\bz)d\bz.
\end{align*}

The fully discretized scheme for the system with random inputs is similar to how we solve the deterministic problem, 
introduced in section \ref{Det_scheme}, except that each term now is a vector analogy of the corresponding term in the deterministic problem. 

\subsection{A Stochastic AP Time-splitting}
\label{sec:sAP}
{\color{blue}Jin, Xiu and Zhu first introduced the notion of stochastic AP (s-AP) in the SG setting \cite{XZJ}. }
s-AP schemes require that as $\varepsilon\to 0$, the SG for the model with uncertainties 
($\mathcal F_z^{\varepsilon}$) automatically becomes a SG approximation for the limiting stochastic diffusion equation 
($\mathcal F_z^0${\color{red}), which is the bipolar drift-diffusion equations (\ref{diff1}) in our case. }
In this section, we formally prove that the time-splitting scheme (\ref{s_relax1})--(\ref{s_relax2}) and 
(\ref{s_trans1})--(\ref{s_trans2}) satisfies the s-AP property. 

As $\epsilon\to 0$, (\ref{s_relax1}) becomes 
\begin{equation}\label{s1}
{\bf \hat r}_{i}= {\boldsymbol{\hat\rho}}_{i}M_i\,, 
\end{equation}
a result proved in Lemma $3$ of \cite{JinLiu}. 
The second equation (\ref{s_relax2}) gives 
\begin{equation}\label{s2}
(\hat j_i)_k = -\sum_l (H_i^{-1})_{kl} \left[ v\cdot \nabla_x(\hat r_i)_l \mp 
\sum_m\sum_n \hat E_m\nabla_v (\hat {r_i})_n\, G_{mnl}\right]. 
\end{equation}
Inserting (\ref{s1}) and (\ref{s2}) into (\ref{s_trans1}) and integrating over $v\in\mathbb R$, one gets
\begin{equation}
\label{gPC_limit}\partial_t (\hat\rho_i)_k - \nabla_x \cdot\left[T_i \sum_l (H_i^{-1})_{kl} \bigg(\nabla_x (\hat\rho_i)_l \pm 
\sum_m\sum_n \hat E_m  (\hat\rho_i)_n\, G_{mnl}\bigg)\right]
=\int_{\mathbb R} ({\bf I}_{\text{i,plus}})_k\, dv, 
\end{equation}
where $$T_i = \int_{\mathbb R} v\otimes v M_i(v)dv, $$
\begin{align}
&\displaystyle  \int_{\mathbb R} ({\bf I}_{\text{1,plus}})_k\, dv = \frac{1}{2}\int_{I_{\bz}}\int_{\mathbb R}\int_{\mathbb R}M_1(v) 
\big(\sigma_I(x,v,w,\bz) + \sigma_I(x,-v,w,\bz)\big)dwdv \,\psi_k(\bz)\pi(\bz)d\bz \notag\\[6pt]
&\label{int_I1}\displaystyle \qquad\qquad\qquad\qquad-\sum_m \sum_n (\hat\rho_1)_m (\hat\rho_2)_n \int_{\mathbb R}\int_{\mathbb R} M_1(v)M_2^2(w) F_{mnk}(x,v,w)\, dwdv, \\[6pt]
&\displaystyle  \int_{\mathbb R} ({\bf I}_{\text{2,plus}})_k\, dv = \frac{1}{2}\int_{I_{\bz}}\int_{\mathbb R}\int_{\mathbb R}M_1(w) 
\big(\sigma_I(x,v,w,\bz) + \sigma_I(x,-v,w,\bz)\big) dwdv \,\psi_k(\bz)\pi(\bz)d\bz \notag \\[6pt]
&\label{int_I2}\displaystyle \qquad\qquad\qquad\qquad -\sum_m \sum_n (\hat\rho_1)_m (\hat\rho_2)_n \int_{\mathbb R}\int_{\mathbb R} M_1(w)M_2^2(v) F_{mnk}(x,v,w)\, dwdv. 
\end{align}
It is obvious that $\int_{\mathbb R}\, ({\bf I}_{\text{1,plus}})_k\, dv= \int_{\mathbb R}\, ({\bf I}_{\text{2,plus}})_k\, dv$. 

On the other hand, applying the ansatz 
$$\rho_i^K =\sum_{k=1}^K (\hat{\rho_i})_k \psi_k(\bz) = {\boldsymbol{\hat\rho}}_{i}\cdot\boldsymbol{\psi}(\bz),  \qquad 
E^K =\sum_{k=1}^K \hat{E}_k \psi_k(\bz), $$
and conducting the Galerkin projection for the limiting drift-diffusion system (\ref{diff1}), one obtains
\begin{equation}
\label{diff_Gal}
\partial_t (\hat\rho_i)_k - \nabla_x \cdot \left[ T_i \sum_{l} (S_i)_{kl}\bigg(\nabla_{x}(\hat\rho_i)_l \pm \sum_m\sum_n \hat E_m (\hat\rho_i)_n\, G_{mnl}\bigg)
\right] = {\bf R}_{k}({\boldsymbol{\hat\rho}}_1, {\boldsymbol{\hat\rho}}_2),  
\end{equation}
where 
$$(S_i)_{kl}=\int_{I_{\bz}}\frac{1}{\lambda_i(x,v,\bz)}\psi_k(\bz)\psi_l(\bz)\pi(\bz)d\bz,$$ 
with $\lambda_i$ defined in (\ref{freq1}), and 
\begin{align}
&\displaystyle  {\bf R}_{k}({\boldsymbol{\hat\rho}}_1, {\boldsymbol{\hat\rho}}_2) = \int_{I_{\bz}}\int_{\mathbb R}\int_{\mathbb R}
\sigma_I(v,w,\bz)M_1(v)dvdw\,\psi_k(\bz)\pi(\bz)d\bz \notag\\[4pt]
&\label{vec_R}\displaystyle \qquad\qquad\qquad - \sum_m\sum_n (\hat\rho_1)_m (\hat\rho_2)_n \int_{\mathbb R}\int_{\mathbb R} 
M_1(v)M_2^2(w) F_{mnk}(x,v,w)dvdw. 
\end{align}
By a change of variable $w'=-w$, the first terms of (\ref{int_I1}), (\ref{int_I2}) and (\ref{vec_R}) are all equal, 
\begin{align*}
&\quad\displaystyle\frac{1}{2}\int_{I_{\bz}}\int_{\mathbb R}\int_{\mathbb R}M_1(v) 
\big(\sigma_I(x,v,w,\bz) + \sigma_I(x,-v,w,\bz)\big)dwdv \,\psi_k(\bz)\pi(\bz)d\bz \\[4pt]
&\displaystyle=\int_{I_{\bz}}\int_{\mathbb R}\int_{\mathbb R}\sigma_I(x,v,w,\bz)M_1(v)dvdw \, \psi_k(\bz)\pi(\bz)d\bz, 
\end{align*}
thus the right-hand-side of (\ref{gPC_limit}) and (\ref{diff_Gal}) are equal. 
We observe that the limiting scheme of gPC-SG method given by (\ref{gPC_limit}) is 
almost exactly the same as the Galerkin system of the bipolar drift-diffusion equations given by (\ref{diff_Gal}), except for the diffusion coefficient 
matrix $(H_i)^{-1}$ and $S_i$. 
It has been demonstrated in \cite{JinLiu} that the matrix $(S_i)_{K\times K} \sim (H_i)^{-1}_{K\times K}$
with spectral accuracy, thus (\ref{gPC_limit}) is a good approximation of (\ref{diff_Gal}). 

This formally shows that with the deterministic AP solver introduced in section \ref{sec:3}, 
the fully discrete time and space approximations of the corresponding gPC-SG scheme introduced in section \ref{subsec:gPC} are s-AP, 
implying that as $\varepsilon\to 0$, with $\Delta t$, $\Delta x$ fixed, 
the gPC-SG scheme approaches the fully discrete gPC-SG approximation of the bipolar drift-diffusion equations. This will be demonstrated 
in our numerical tests. 

{\color{red} {\bf Remark.} 
With the non-linear generation-recombination integral operators, the proof is 
different from the previous work \cite{JinLiu}, where the gPC-SG scheme for the linear 
semiconductor Boltzmann equation with random inputs is studied. 
}

\section{Numerical examples}
\label{sec:6}

In this section, several numerical tests are shown to illustrate the validity and effectiveness of our AP scheme for the deterministic problem 
(Test 1) and for the model with uncertainties (Test 2). 

In application, people often are more interested in the solution statistics, such as the mean and standard deviation of the macroscopic physical quantities. 
The macroscopic quantities $\rho$, $\mu$ that stand for density and bulk velocity are defined by
\begin{equation}\label{macro} \rho = \int_{\mathbb R} f(v)dv, \qquad \mu=\frac{1}{\rho}\int_{\mathbb R}f(v)v dv,  \end{equation}
and we denote momentum $u=\int_{\mathbb R}f(v)v dv$ in the figures. 

Given the gPC coefficients $f_k$ of $f$, the statistical mean, variance and standard deviation are 
$$ E[f]\approx f_1,  \qquad\text{Var}[f]\approx\sum_{k=2}^{K} f_{k}^2,  \qquad SD[f]=\sqrt{\sum_{k=2}^{K} f_{k}^2}\,. $$

The computational domain is $x\in [0, 1]$ for all the numerical tests. $i=1$ stands for the electrons and $i=2$ stands for the holes. 

\newpage
{\bf{\large Test 1: The deterministic model}}

The equilibrium boundary condition in $x$ is assumed, 
$$f_i(x_L, v, t)=M_i(v),\, v>0\, ;  \qquad f_i(x_R, v, t)=M_i(v), \,  v<0\,. $$
The initial distribution is $f_i(x,v,t=0)=M_i(v)$, for $i=1, 2$. 
The collision and generation-recombination kernels are given by
$$\sigma_1(v,w)=\sigma_2(v,w)=2,  \qquad \sigma_{I}(v,w)=\frac{1}{\sqrt{\pi}}e^{-(v-w)^2}, $$
and 
$$ \beta=0.9, \qquad \gamma=0.002, \qquad \Phi(x_L)=0, \qquad \Phi(x_R)=5, $$
where $\Phi(x_L), \Phi(x_R)$ are the boundary data of the potential at $x_L$, $x_R$ respectively. 
$$ c(x) = 1-(1-m)\left[\tanh(\frac{x-x_1}{s})-\tanh(\frac{x-x_2}{s})\right], $$
with $s=0.02$, $m=(1-0.001)/2$, $x_1=0.3$, $x_2=0.7$. The parameters are chosen similarly as \cite{JinLorenzo}. 
\\[20pt]
{\bf\large{Test 1 a): Convergence to the equilibrium test}}

Denote the discretized numerical solution $f_i(x_l, v_m, T)$ and $\rho_i(x_l, T)$ by $f_i^{l,m}$ and $\rho_i^{l}$ ($i=1,2$), 
where $0\leq l\leq N_x$, $0\leq m\leq N_v$, with $N_x$, $N_v$ the number of mesh points used in $x$ and $v$ directions respectively, and 
$T$ is the final computation time. 

\begin{figure}[H]
\centering
 \includegraphics[width=0.496\linewidth]{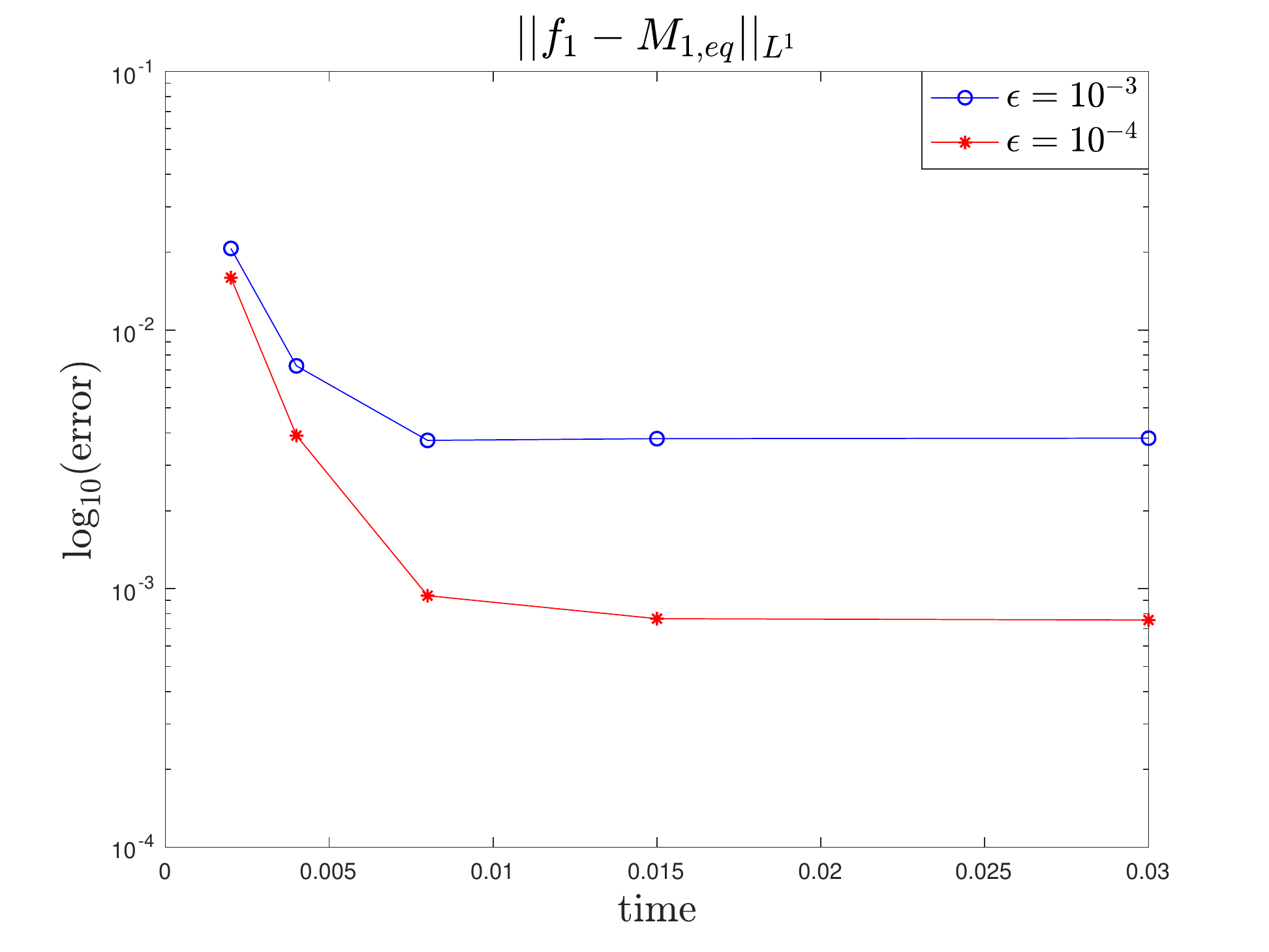}
 \centering
 \includegraphics[width=0.496\linewidth]{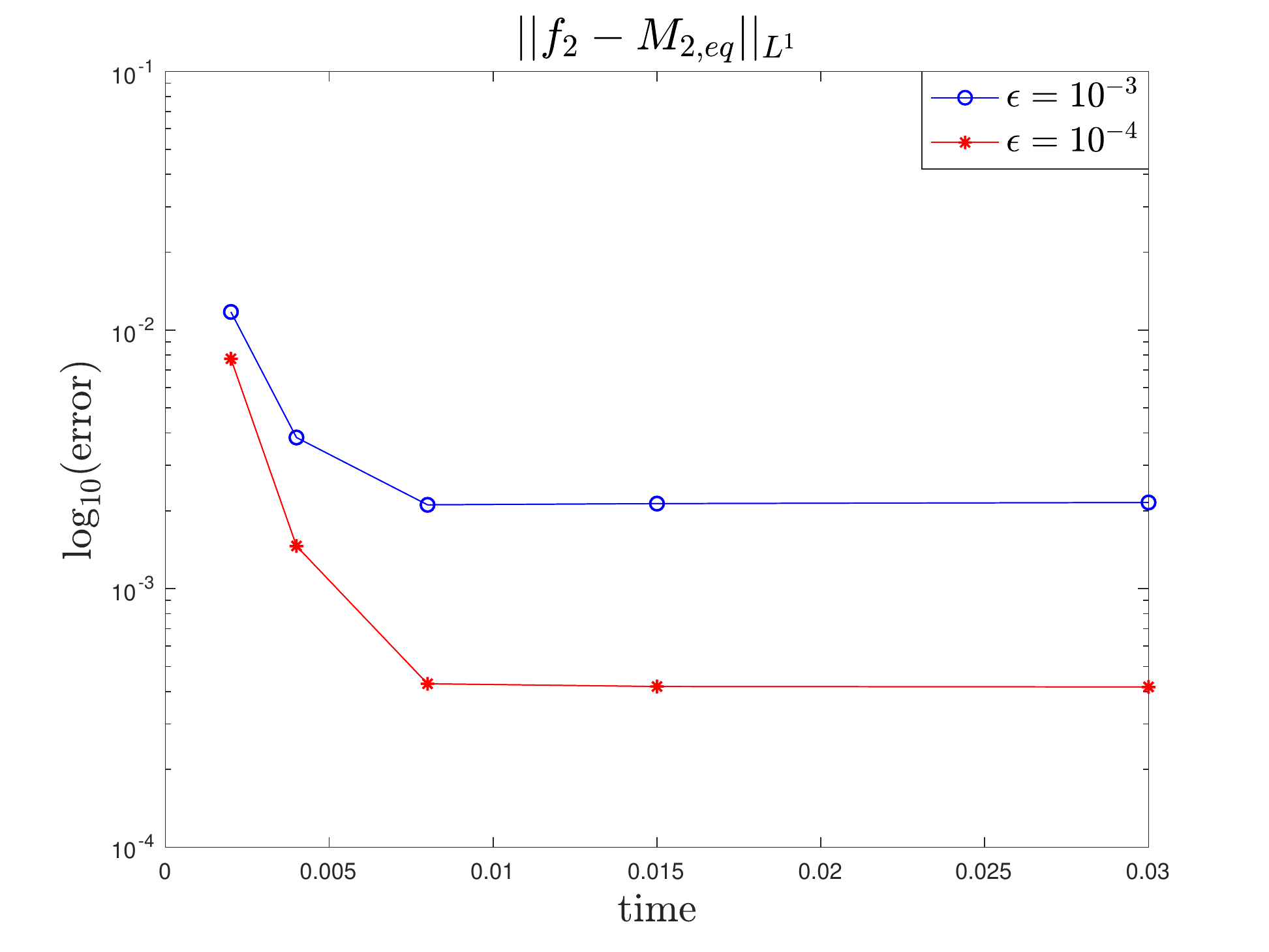}
\caption{Test 1 a). The time evolution of $||f_i -M_{i, \text{eq}}||_{L^1}$ with respect to different $\varepsilon$. 
$\Delta x=0.01$, $N_v=20$, $\Delta t=2\times 10^{-6}$. }
\label{Error_AP}
\end{figure}

Figure \ref{Error_AP} shows the asymptotic error in $L^1(x,v)$ norm by the distance between the distribution function $f_i$ and its corresponding local equilibrium $M_{i, \text{eq}}$ ($i=1, 2$), defined by 
$$ ||f_i - M_{i, \text{eq}}||_{L^1}=||f_i - \rho_i M_i ||_{L^1} = \int_{\mathbb R}\int_{\mathbb R} |f_i-\rho_i M_i|\, dxdv 
=\sum_{l,m}|f_{i}^{l,m} -\rho_i^{l}M_i|\, \Delta x \Delta v, $$
with $M_i$ the absolute Maxwellian given in (\ref{Max}). We report the results for $\varepsilon=10^{-3}$ and $\varepsilon=10^{-4}$. 
As expected, the asymptotic error is $O(\varepsilon)$ before it saturates and the numerical errors
from spacial, temporal and velocity discretizations start to dominate. 

\newpage{\bf\large{Test 1 b): The AP property}}

\begin{figure}[H]
\centering
\includegraphics[width=0.496\linewidth]{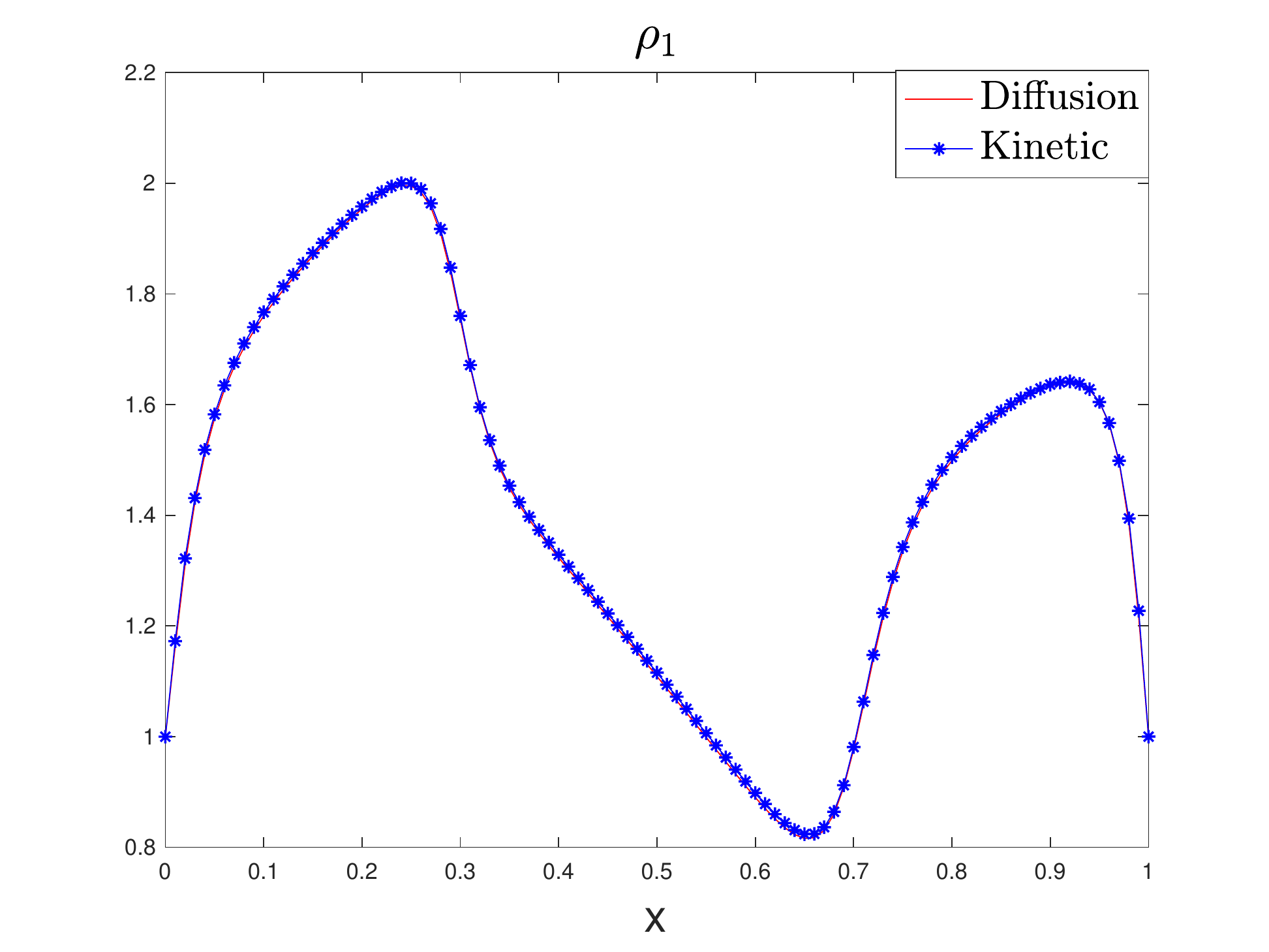}       
\centering
\includegraphics[width=0.496\linewidth]{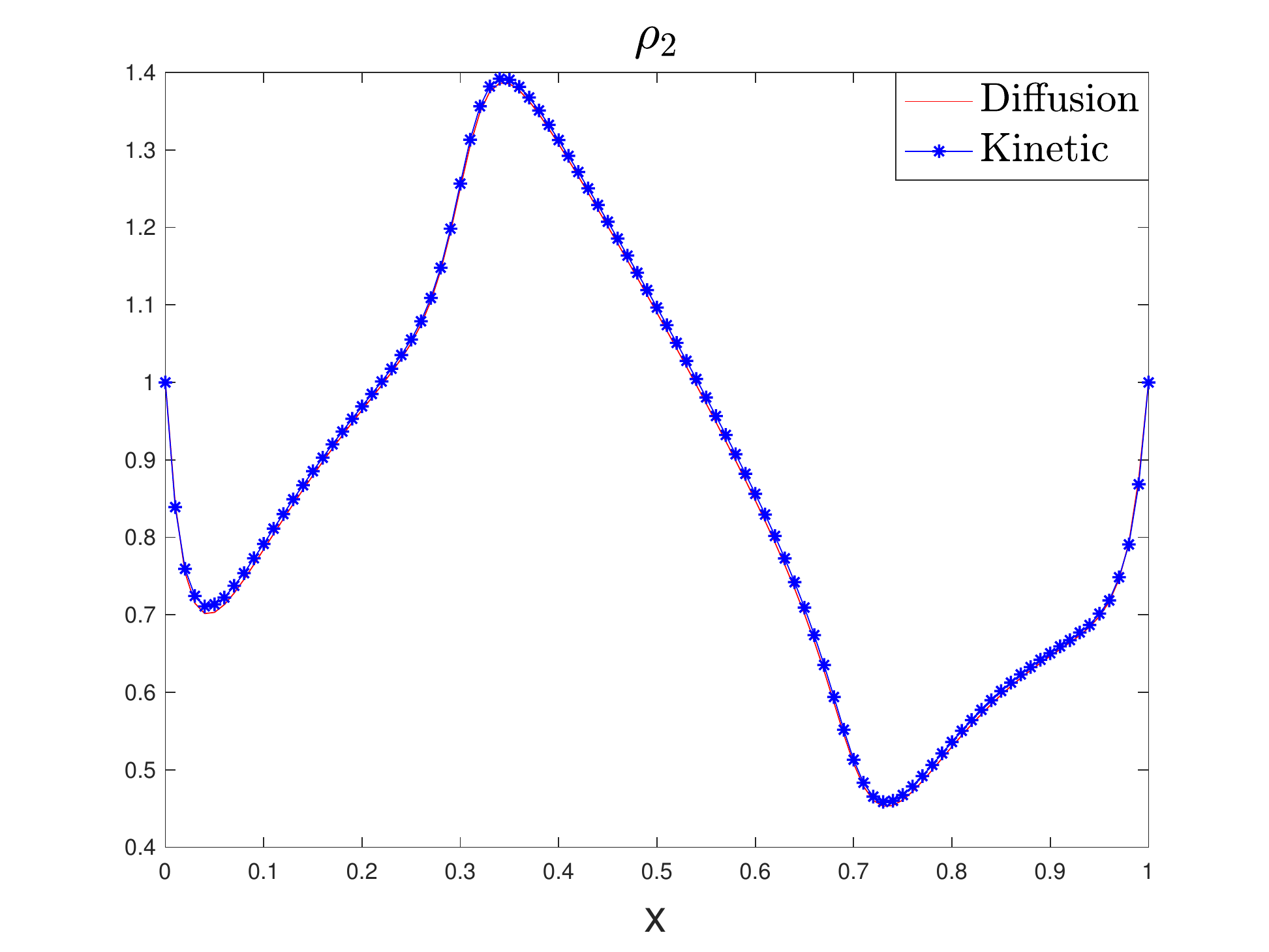}        
\caption{Test 1 b). Solutions at $T=0.2$. $\varepsilon=10^{-5}$, $\Delta x=0.01$, $N_v=20$, $\Delta t=2\times 10^{-6}$ for the kinetic model; 
and $N_v=32$, $\Delta x=5\times 10^{-3}$, $\Delta t=2\times 10^{-6}$ for the drift-diffusion system. }
\label{Test1b}
\end{figure}

Figure \ref{Test1b} demonstrates that when $\varepsilon$ is really small ($\varepsilon=10^{-5}$), 
the solutions of the kinetic system $\rho_1$, $\rho_2$ automatically becomes the solutions of the bipolar drift-diffusion system, 
known as the desired AP property. The forward Euler in time and the central difference scheme in space is used to compute the numerical approximations 
(with fine grids) of the drift-diffusion equations. One can observe that two sets of solutions are in good agreement. 

Test 2 below studies the model with random inputs and validate the efficiency and accuracy of our s-AP gPC-SG method. 

The stochastic collocation (SC) method \cite{XiuBook} is employed for numerical comparison with the gPC-SG method. 
We explain the basic idea. Let $\{ \bz^{(j)}\}_{j=1}^{N_c} \subset I_{\bz}$ be the set of collocation nodes and $N_c$ the number of collocation points. 
For each individual sample $\bz^{(j)}$, $j=1, \cdots, N_c$, one applies the deterministic AP solver to obtain the solution at sampling points
$f_j(t,x,v)=f(t,x,v,\bz^{(j)})$, then adopts the interpolation approach to construct a gPC approximation, such as 
$$ f(t,x,v,\bz)=\sum_{j=1}^{N_c} f_j(t,x,v) l_j(\bz), $$
where $l_j(\bz)$ depends on the construction method. The Lagrange interpolation is used here by choosing $l_j(\bz^{(i)})=\delta_{ij}$. 
In the collocation method, the integrals are approximated by 
$$\int_{I_{\bz}} f(t,x,v,\bz)\pi(\bz)d\bz\approx\sum_{j=1}^{N_c}f(t,x,v,\bz^{(j)})w^{(j)}, $$
where $\{w^{(j)}\}$ are the weights corresponding to the sample points $\{\bz^{(j)}\}$ ($j=1, \cdots, N_c$) from the quadrature rule. 

To measure the difference in mean and standard deviation of the macroscopic quantities given in (\ref{macro}), 
we use $L^2$ norm in $x$ in Test 2 c), 
\begin{align*}\mathcal E_{\text{mean}}(t)=\left|\left|E[w^h]-E[w]\right|\right|_{L^2}, \\[4pt]
\mathcal E_{\text{std}}(t)=\left|\left|SD[w^h]-SD[w]\right|\right|_{L^2}, \end{align*}
where $w^h$ and $w$ are numerical solutions of gPC-SG method and reference solutions obtained by the collocation method. 

In Test 2 a), b), c), we will assume the random variable $z$ obeys a uniform distribution, defined on $[-1,1]$, so the Legendre gPC polynomial basis is 
used. We put different sources of random inputs including the random doping profile, random collision kernels and random initial data
 in Test 2 a), b), c) respectively. We report the results obtained for $\varepsilon=10^{-3}$ at output time $T=0.1$ in Test 2 a), b), c). 
\\[20pt]
{\bf\large{Test 2 a): Random doping profile}}

We assume a random doping profile 
$$ c(x,z) = \left[1-(1-m)\bigg(\tanh(\frac{x-x_1}{s})-\tanh(\frac{x-x_2}{s})\bigg)\right](1+0.5z), $$
and random collision kernels 
$$\sigma_1=\sigma_2=2+z, \qquad \sigma_{I}(v,w)=\frac{1}{\sqrt{\pi}}e^{-(v-w)^2}. $$
Other parameters, initial and boundary data are chosen the same as Test 1. 
 \\[20pt]
{\bf\large{Test 2 b): Random collision kernels}}

Let $$\sigma_1=\sigma_2=2+0.5z, \qquad \sigma_{I}(v,w)=\frac{1}{\sqrt{\pi}}e^{-(v-w)^2}, $$
and other parameters, initial and boundary data are chosen the same as Test 1. 
\\[20pt]
{\bf\large{Test 2 c): Random initial data}}

Assume an initial data with a smooth, random perturbation around its absolute Maxwellian, 
$$ f_i(x,v,t=0)=\rho(z)M_i(v), \qquad \rho(z)=\sin\left[\frac{\pi}{2}(z+1)\right], $$
for $i=1,2$. Other parameters, boundary data are chosen the same as Test 1. 

\begin{figure} [H]
  \centering
 \includegraphics[width=0.49\linewidth]{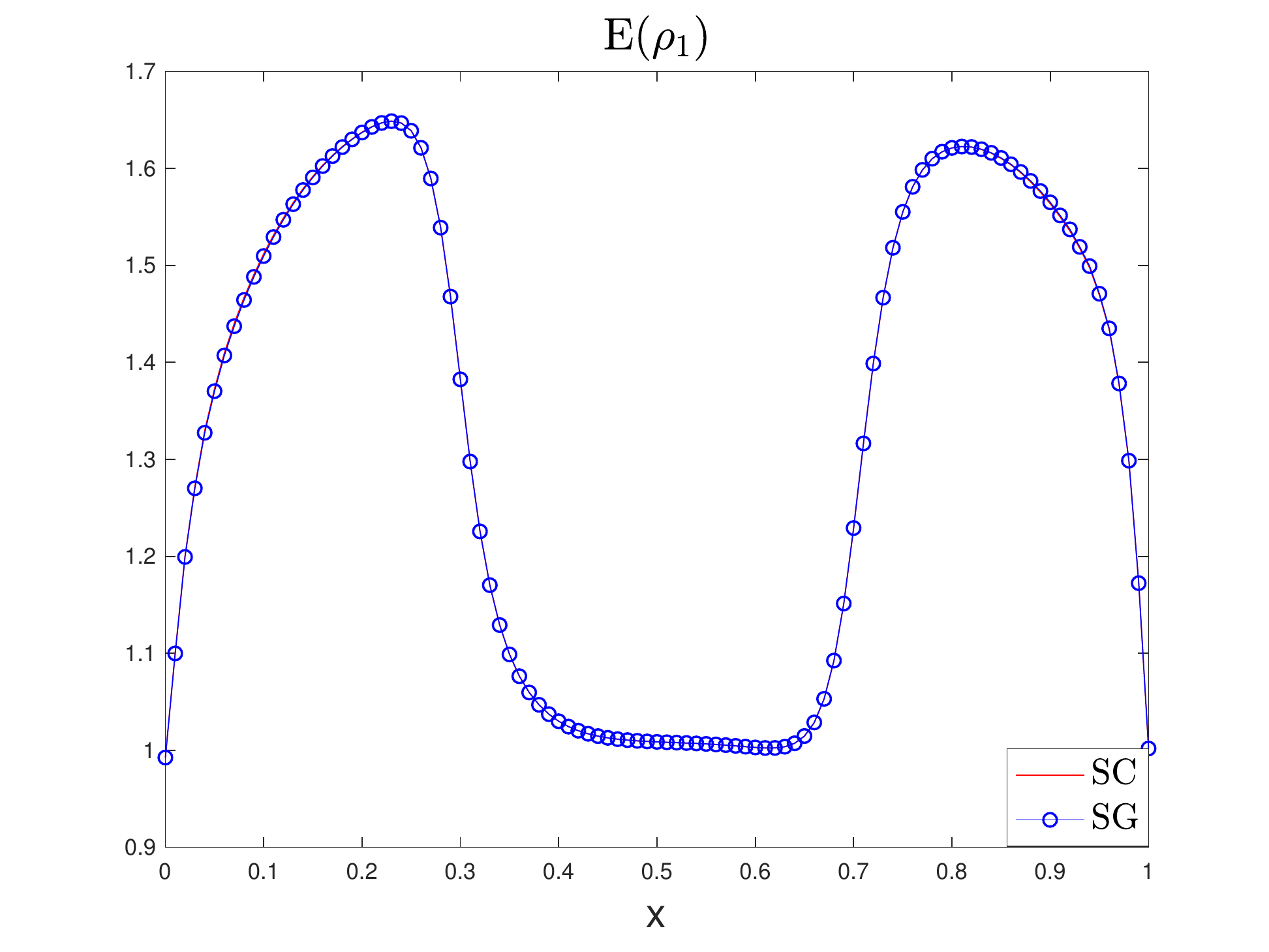}
 \centering
  \includegraphics[width=0.49\linewidth]{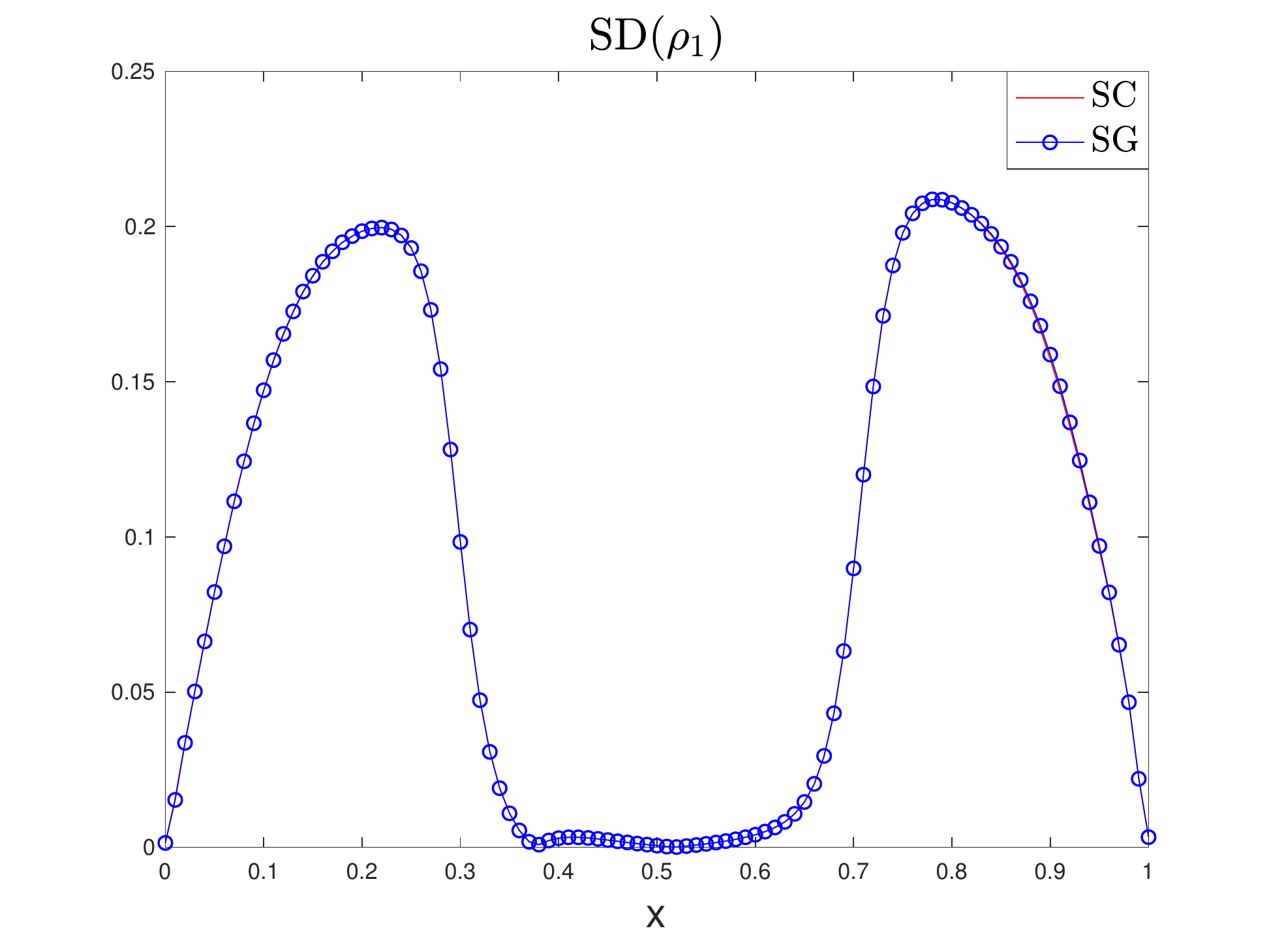}
   \centering
 \includegraphics[width=0.49\linewidth]{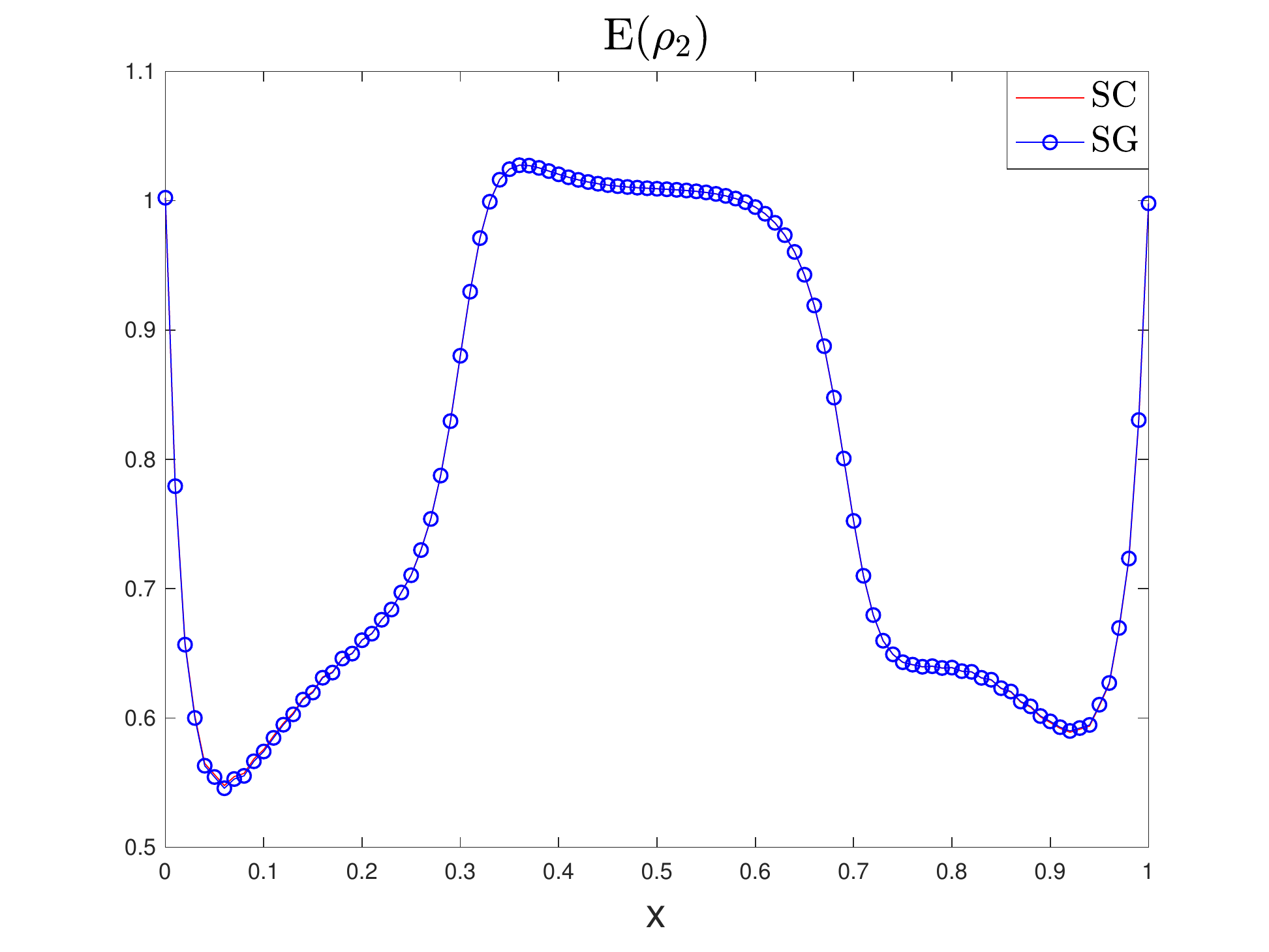}
\centering
\includegraphics[width=0.49\linewidth]{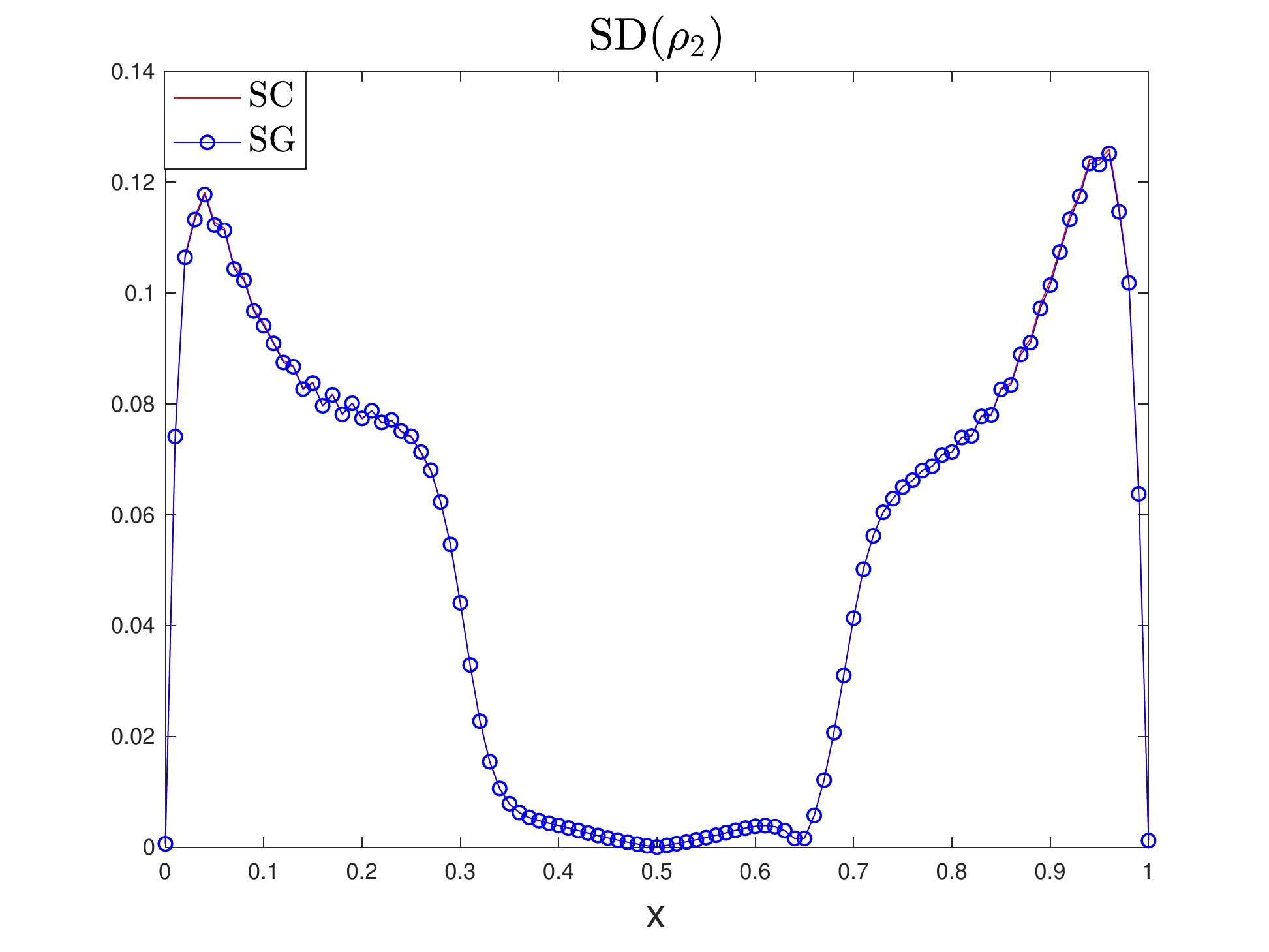}
\centering
 \includegraphics[width=0.49\linewidth]{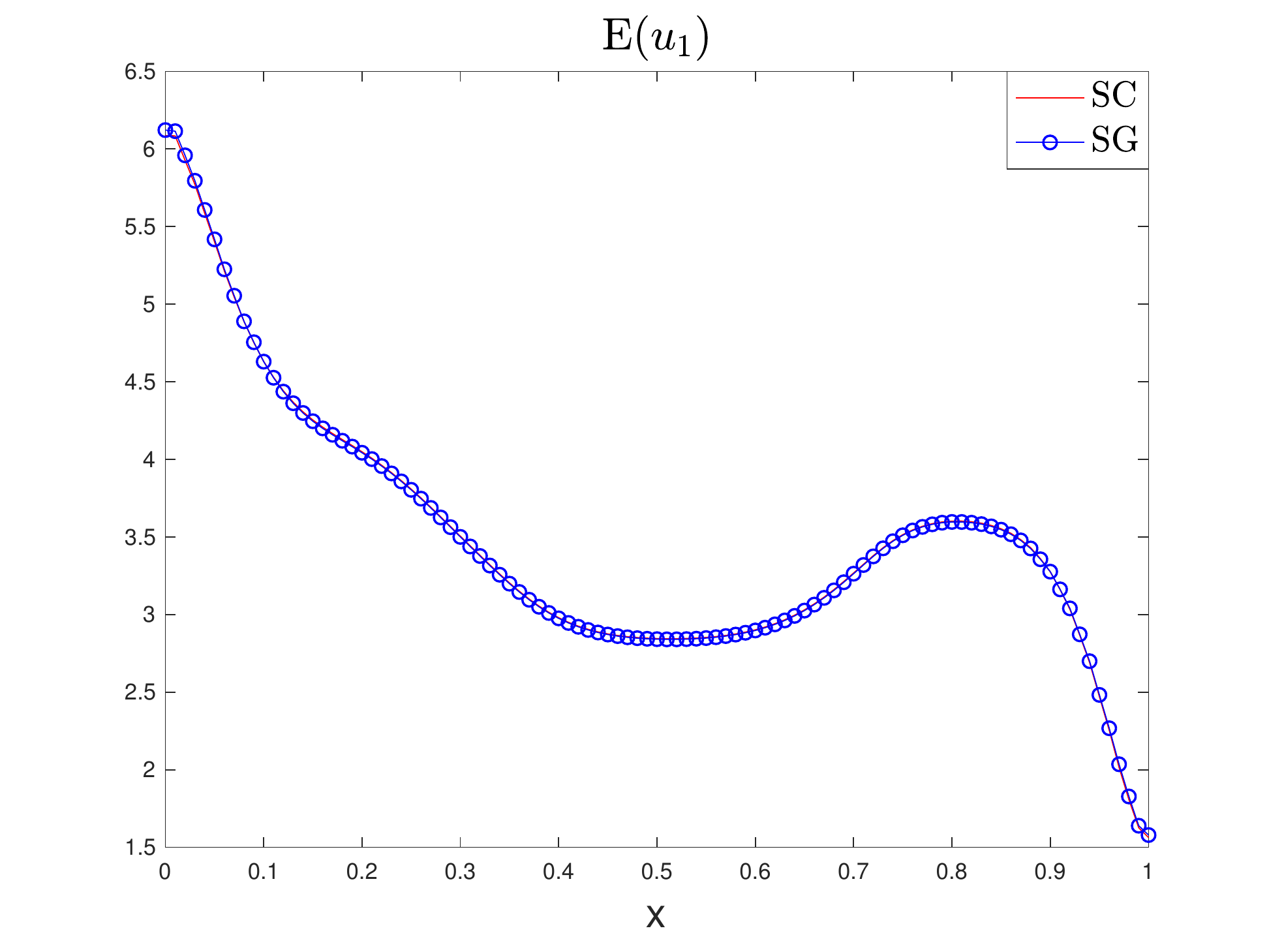}
 \centering
 \includegraphics[width=0.49\linewidth]{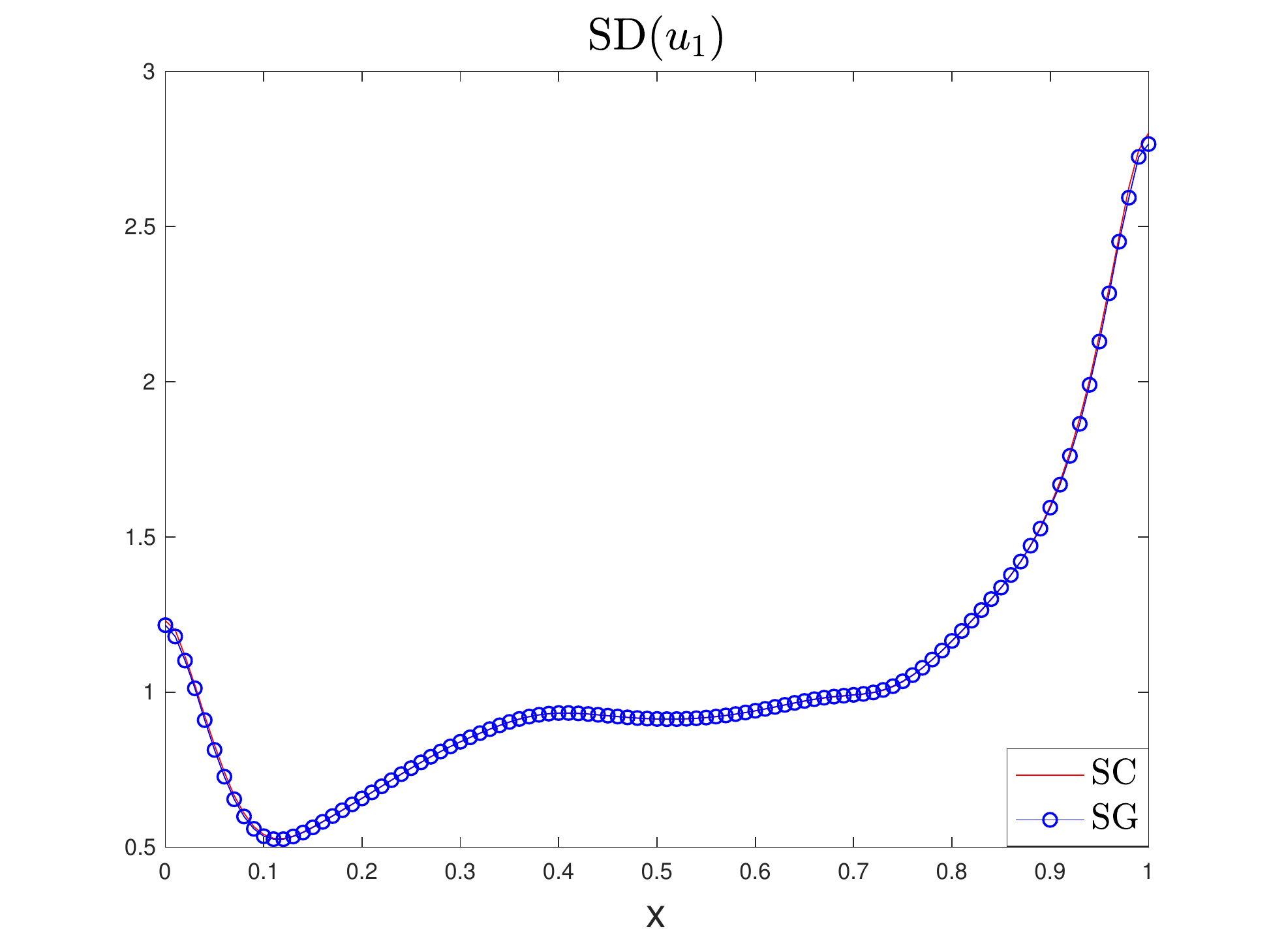}
 \centering
 \includegraphics[width=0.49\linewidth]{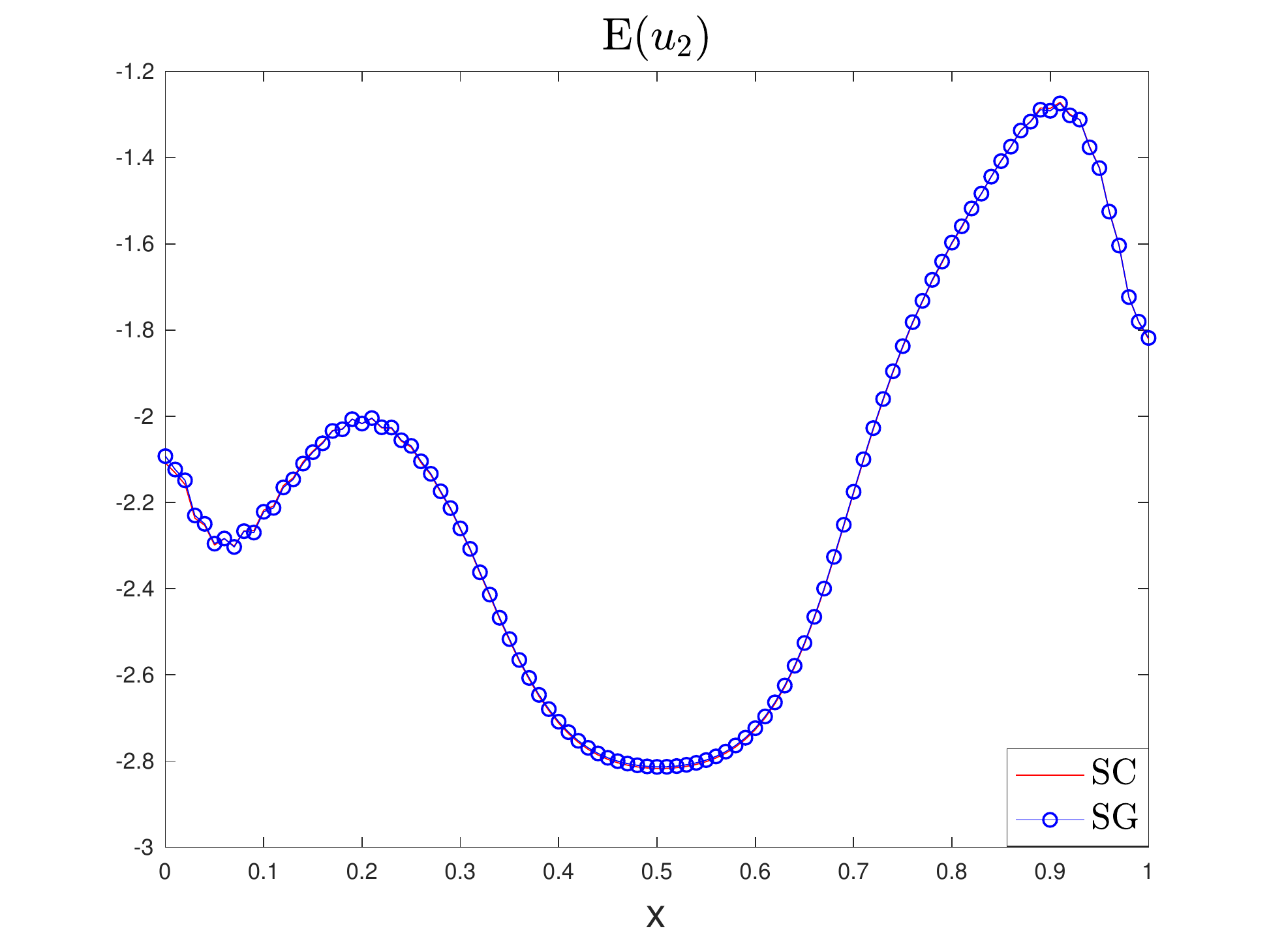}
 \centering
\includegraphics[width=0.49\linewidth]{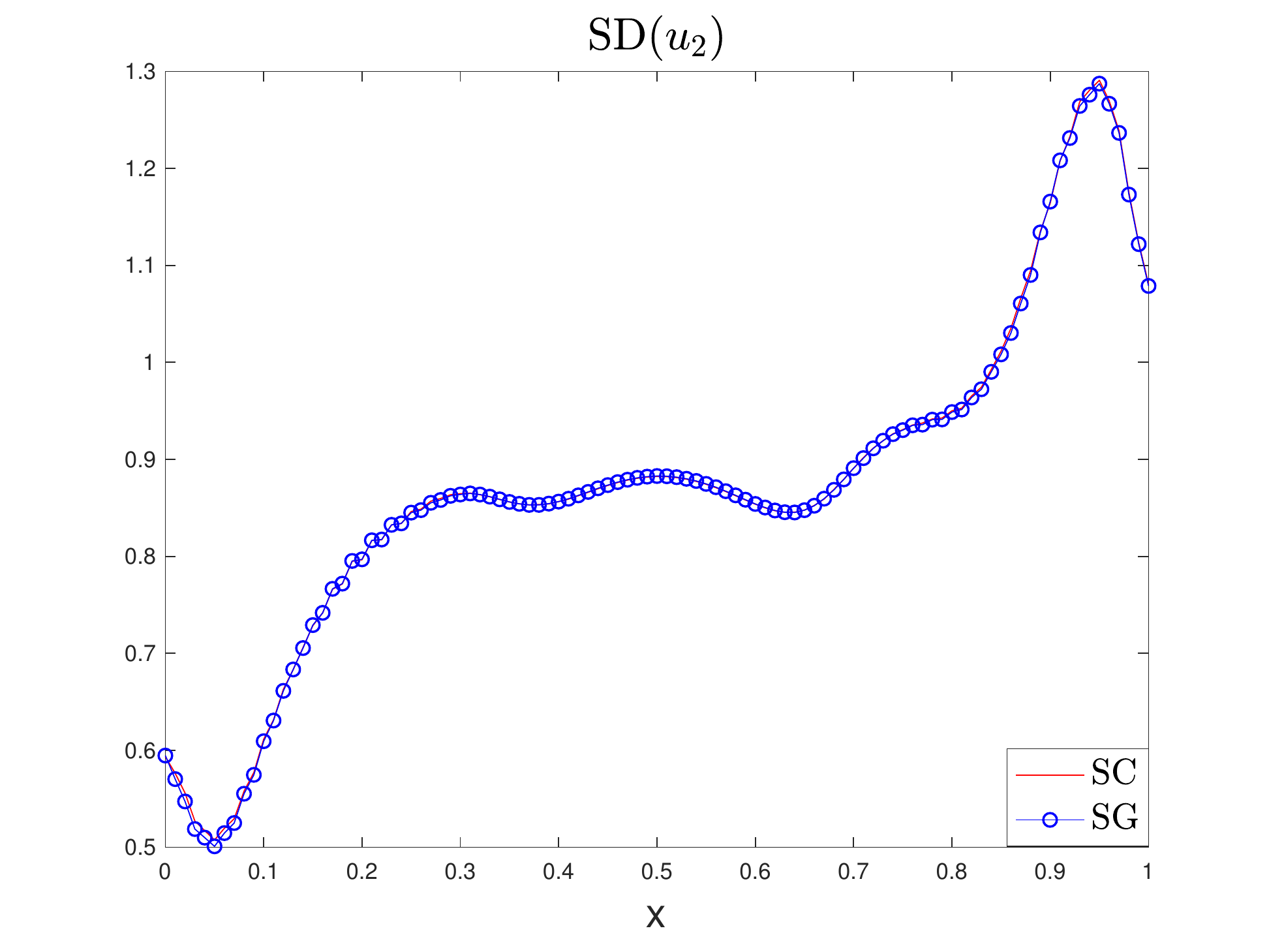}
\caption{Test 2 a). Red solid line: reference solutions by the SC method. Blue line with circles: gPC-SG method with $K=4$. }
\label{Test2a}
\end{figure}

\begin{figure} [H]
  \centering
 \includegraphics[width=0.49\linewidth]{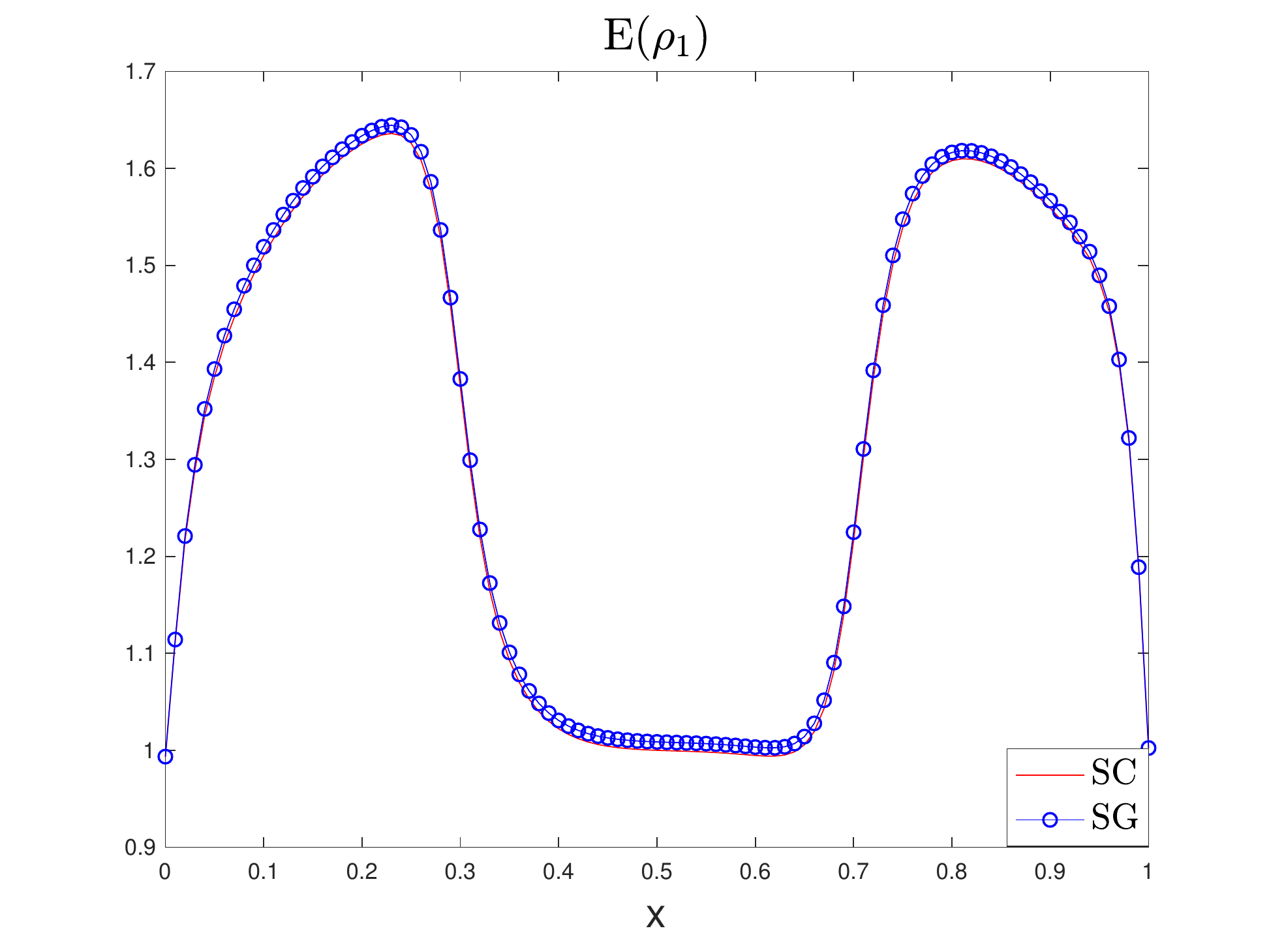}
 \centering
  \includegraphics[width=0.49\linewidth]{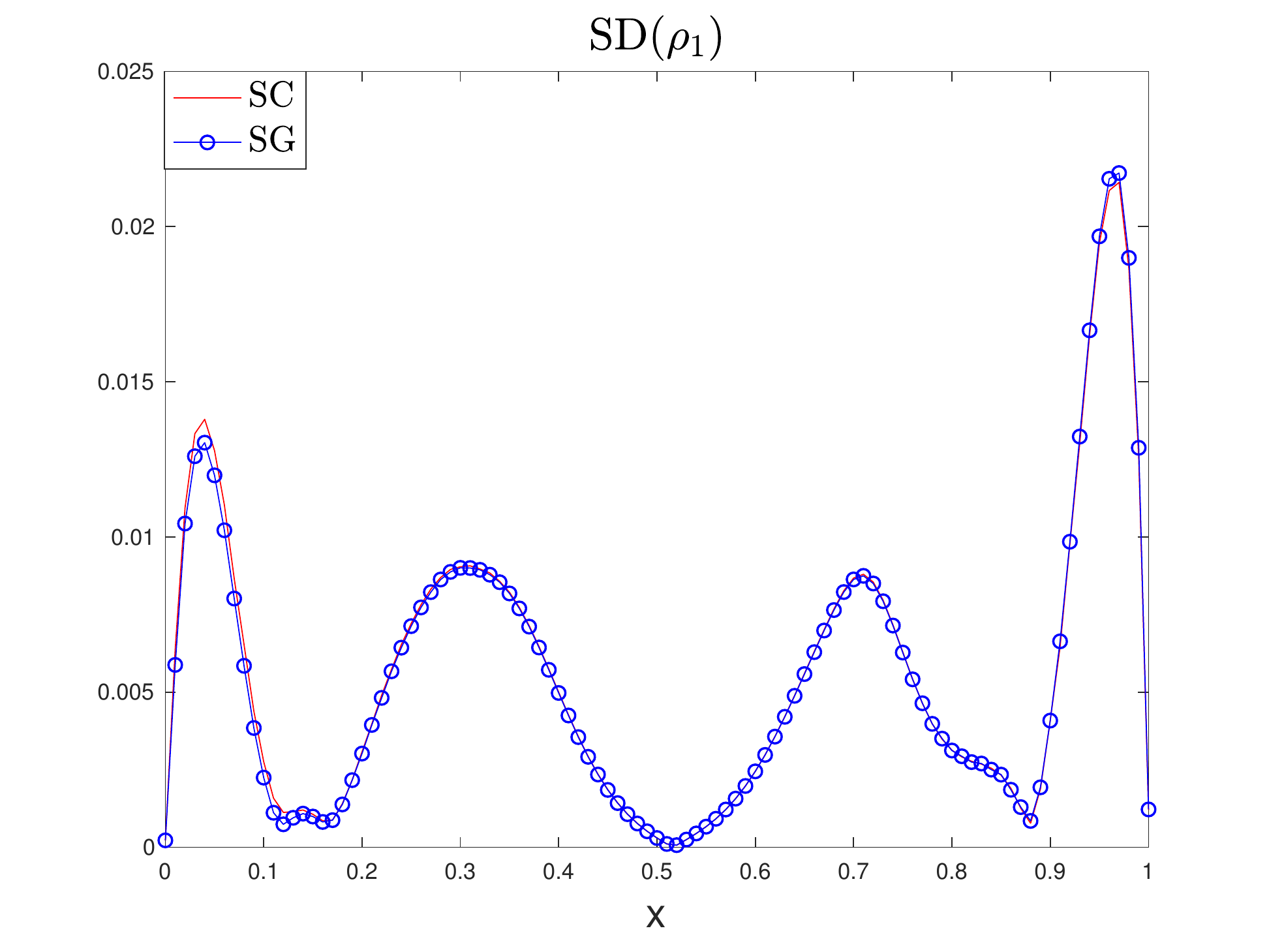}
   \centering
 \includegraphics[width=0.49\linewidth]{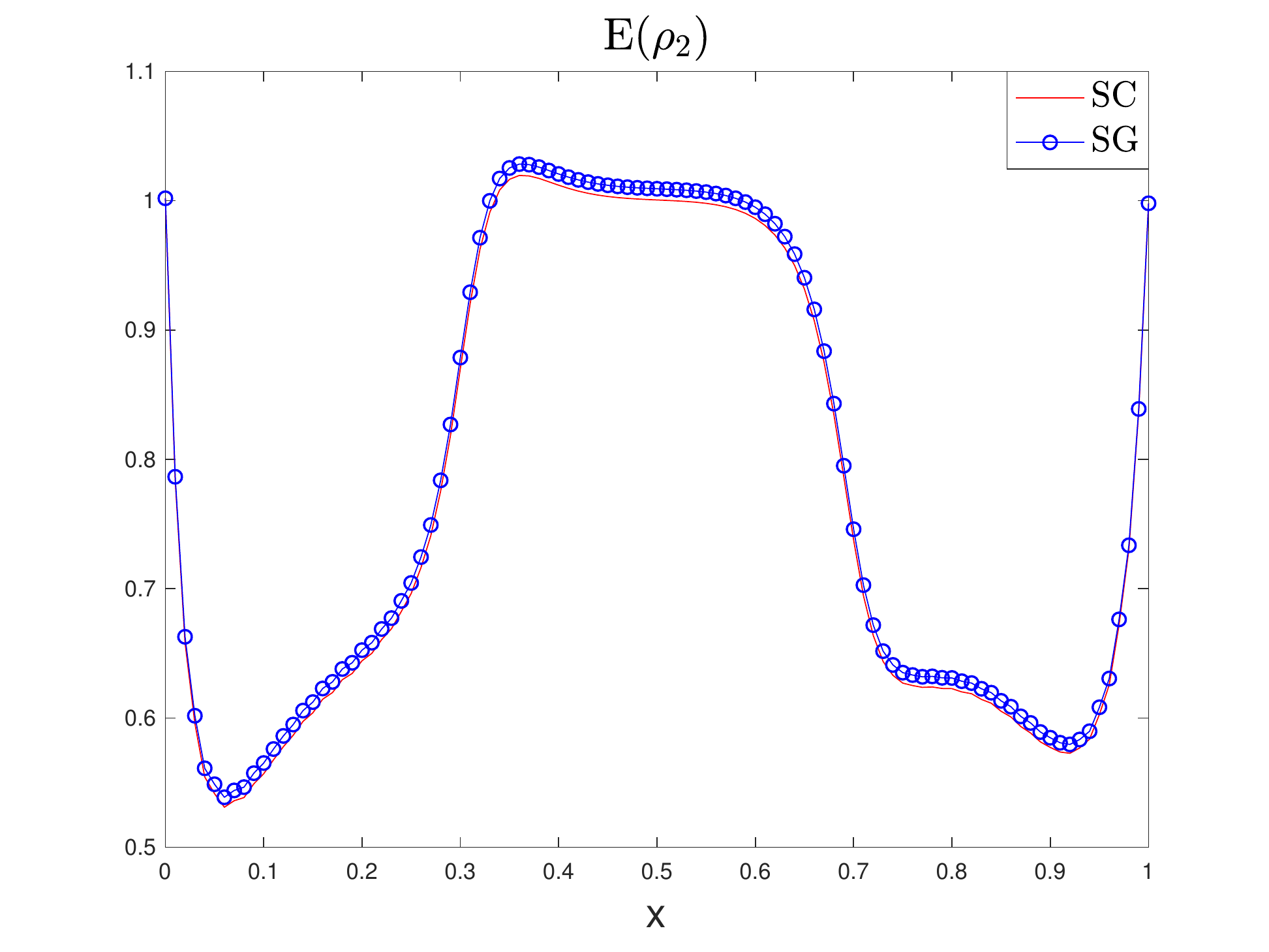}
\centering
\includegraphics[width=0.49\linewidth]{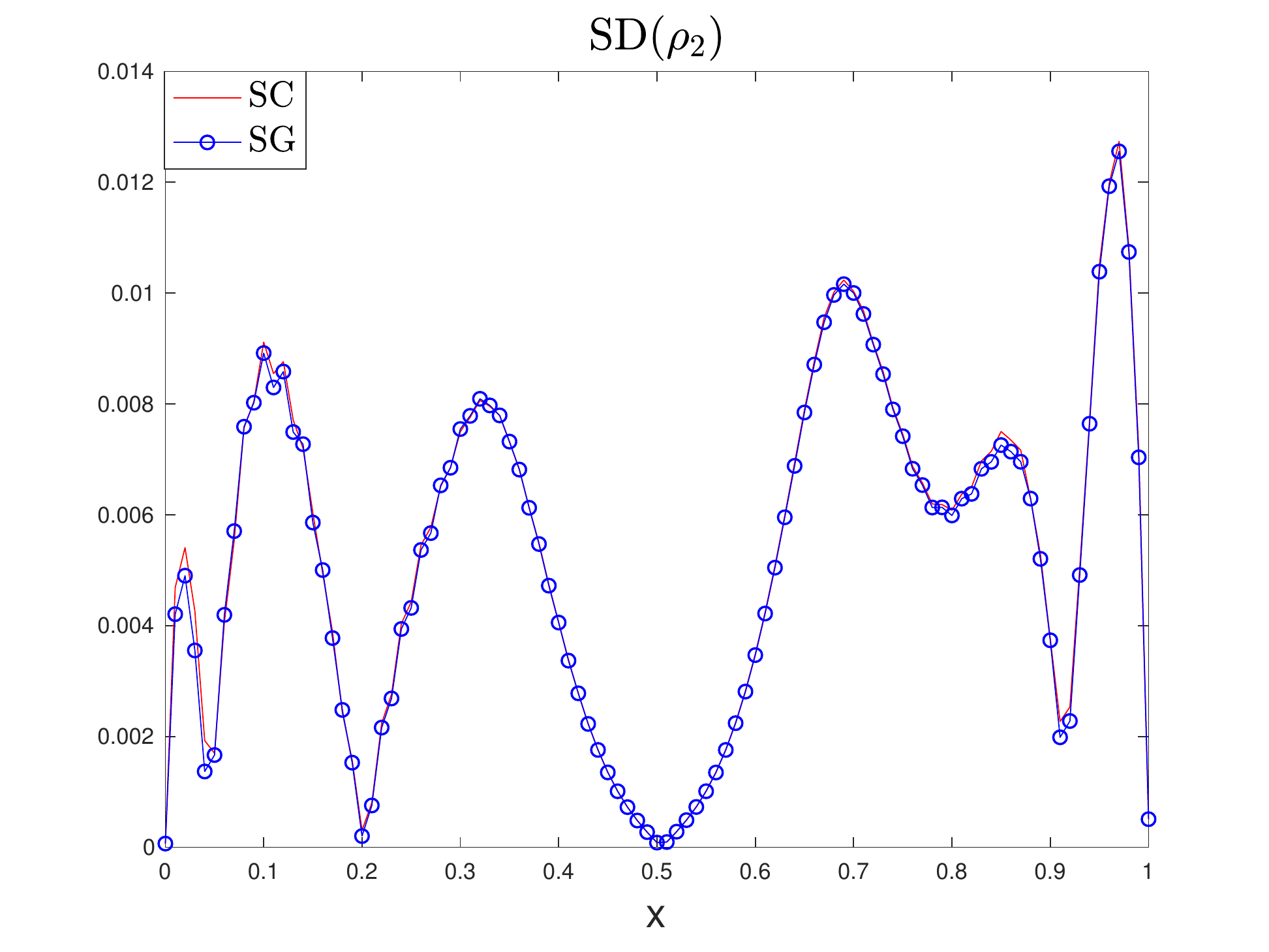}
\centering
 \includegraphics[width=0.49\linewidth]{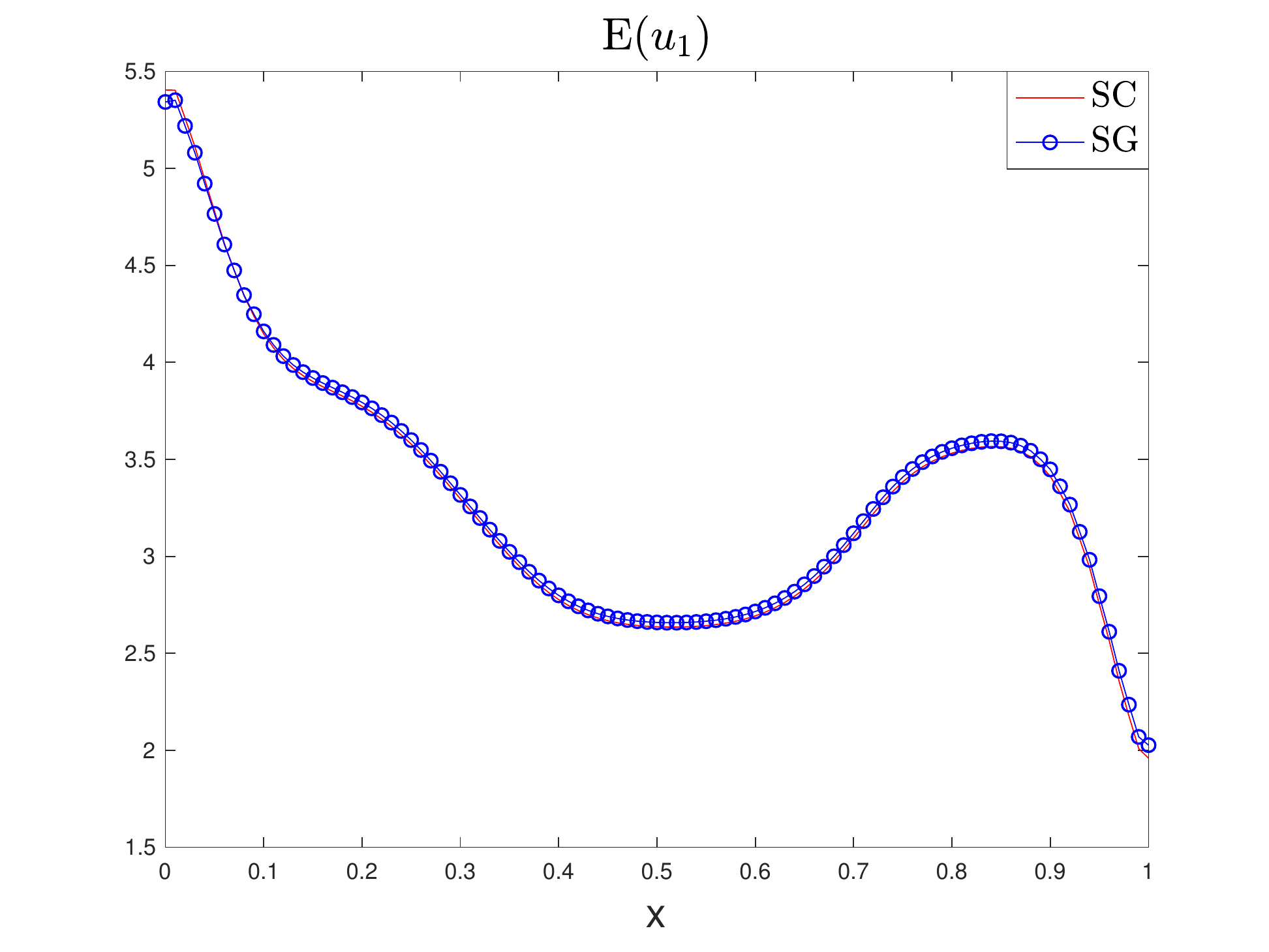}
 \centering
 \includegraphics[width=0.49\linewidth]{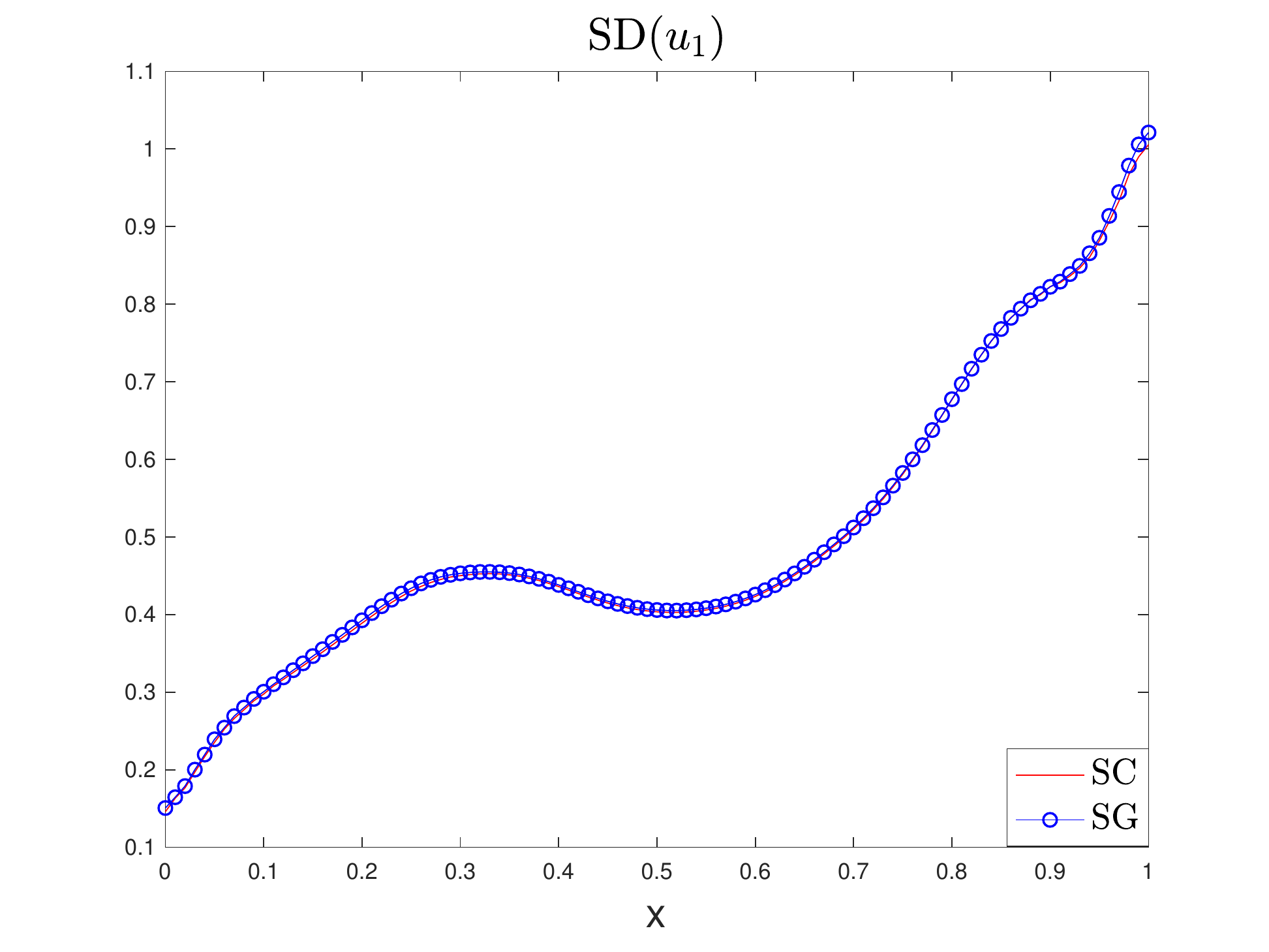}
 \centering
 \includegraphics[width=0.49\linewidth]{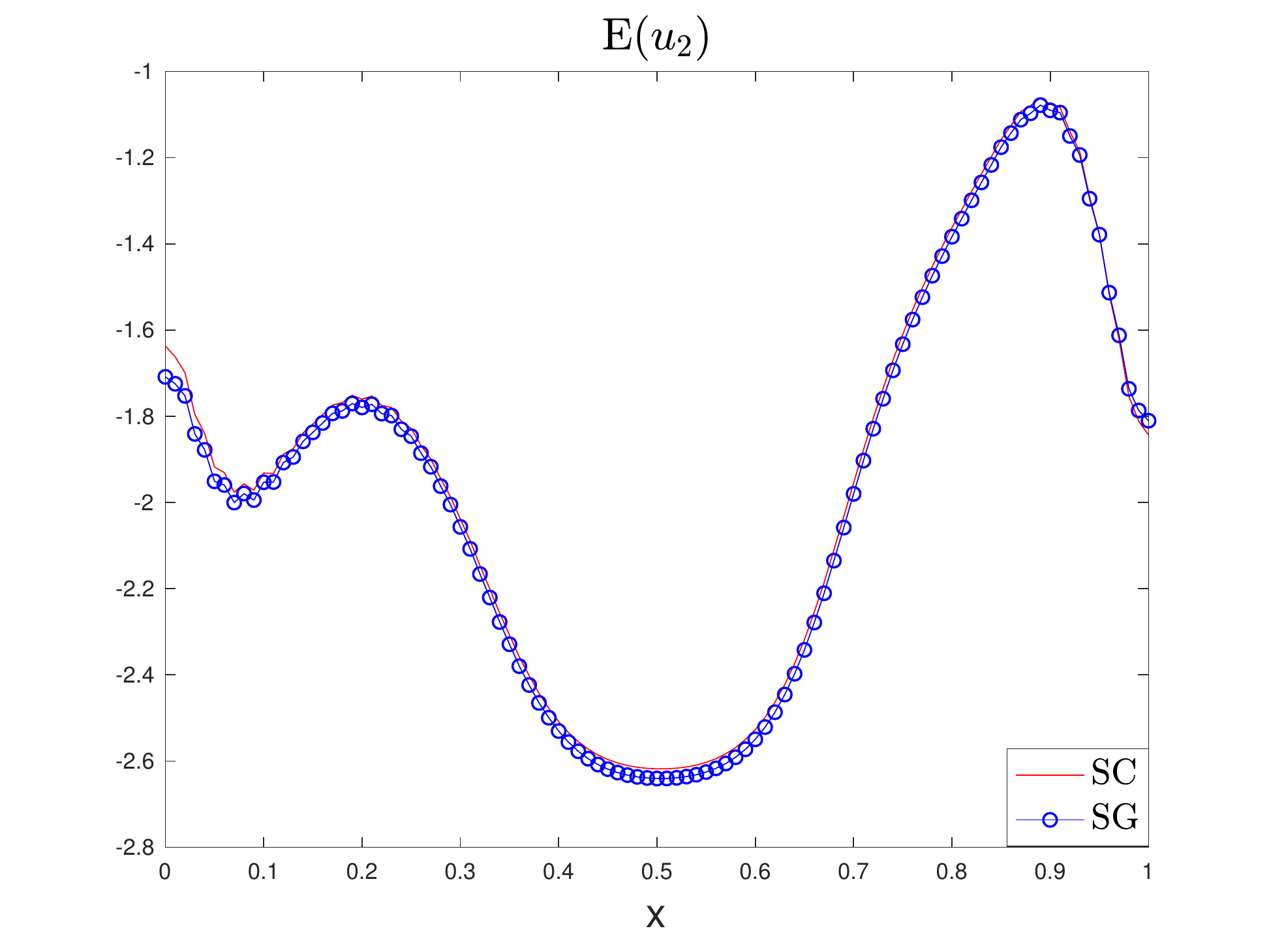}
 \centering
\includegraphics[width=0.49\linewidth]{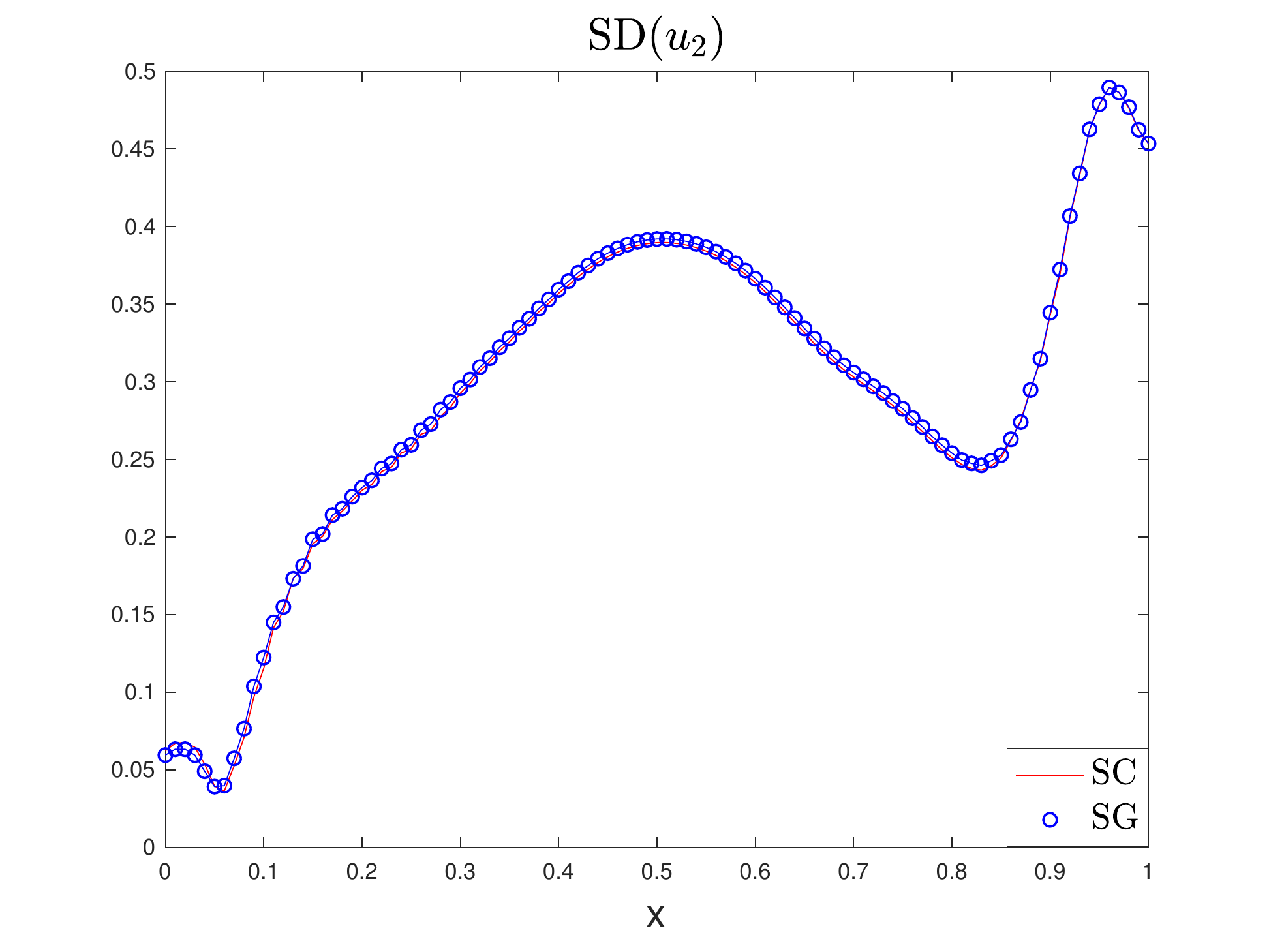}
\caption{Test 2 b).  $\Delta x=0.01, \Delta t=2\times 10^{-6}, N_v=16$. 
Red solid line: reference solutions by the SC method with $N_c=16$.  Blue line with circles: gPC-SG method with $K=4$. }
\label{Test2b}
\end{figure}

\begin{figure} [H]
  \centering
 \includegraphics[width=0.49\linewidth]{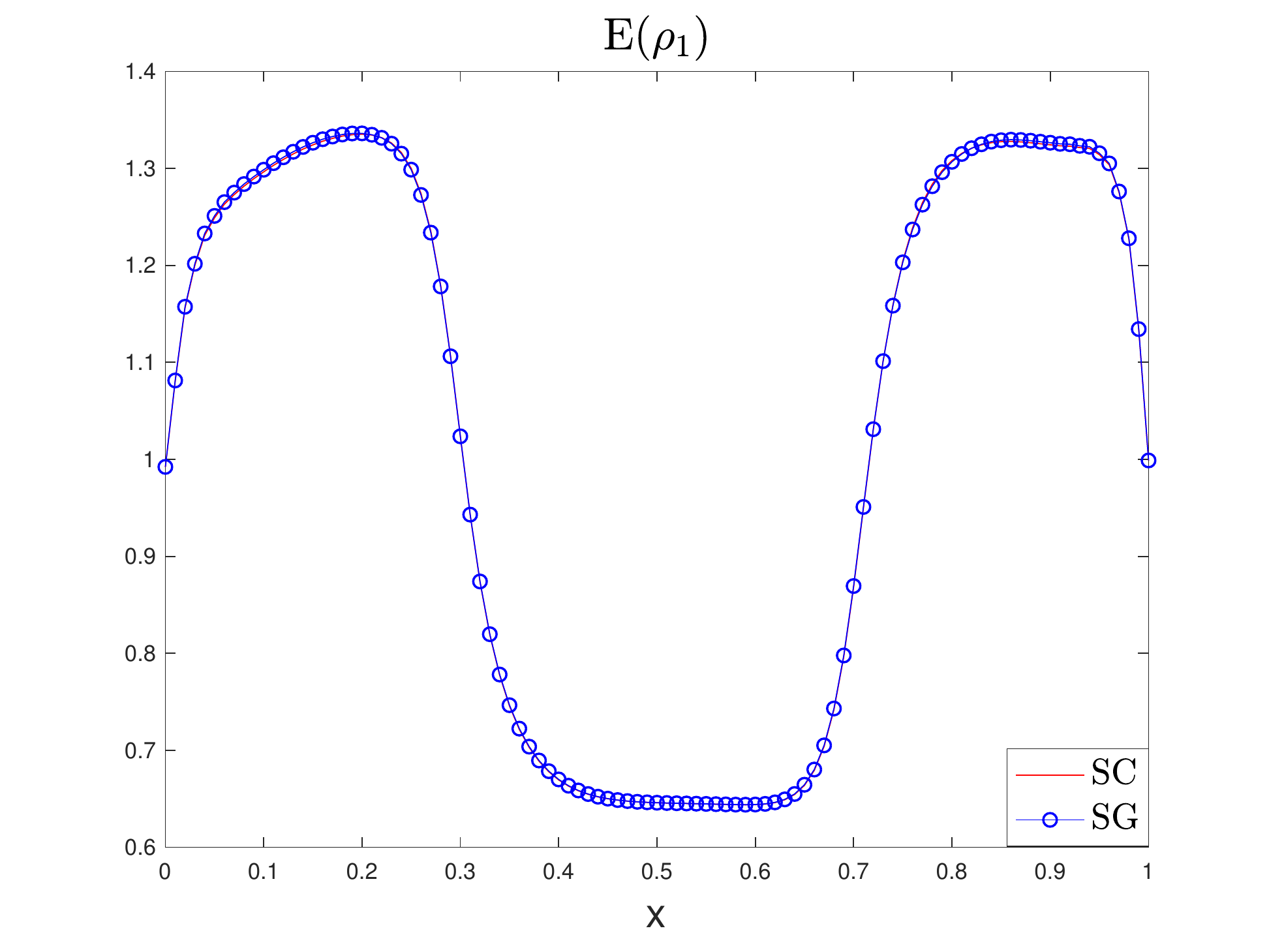}
 \centering
  \includegraphics[width=0.49\linewidth]{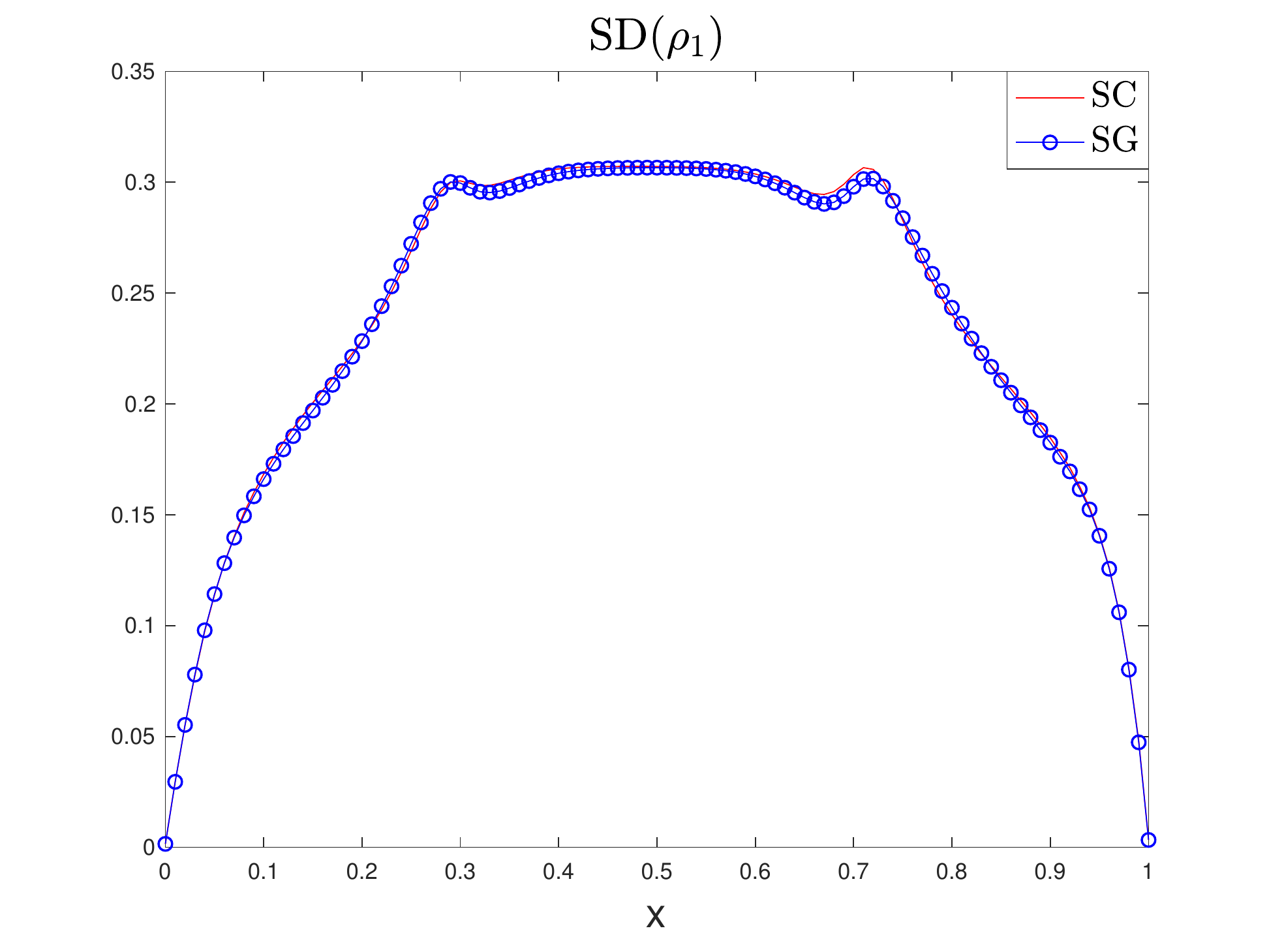}
   \centering
 \includegraphics[width=0.49\linewidth]{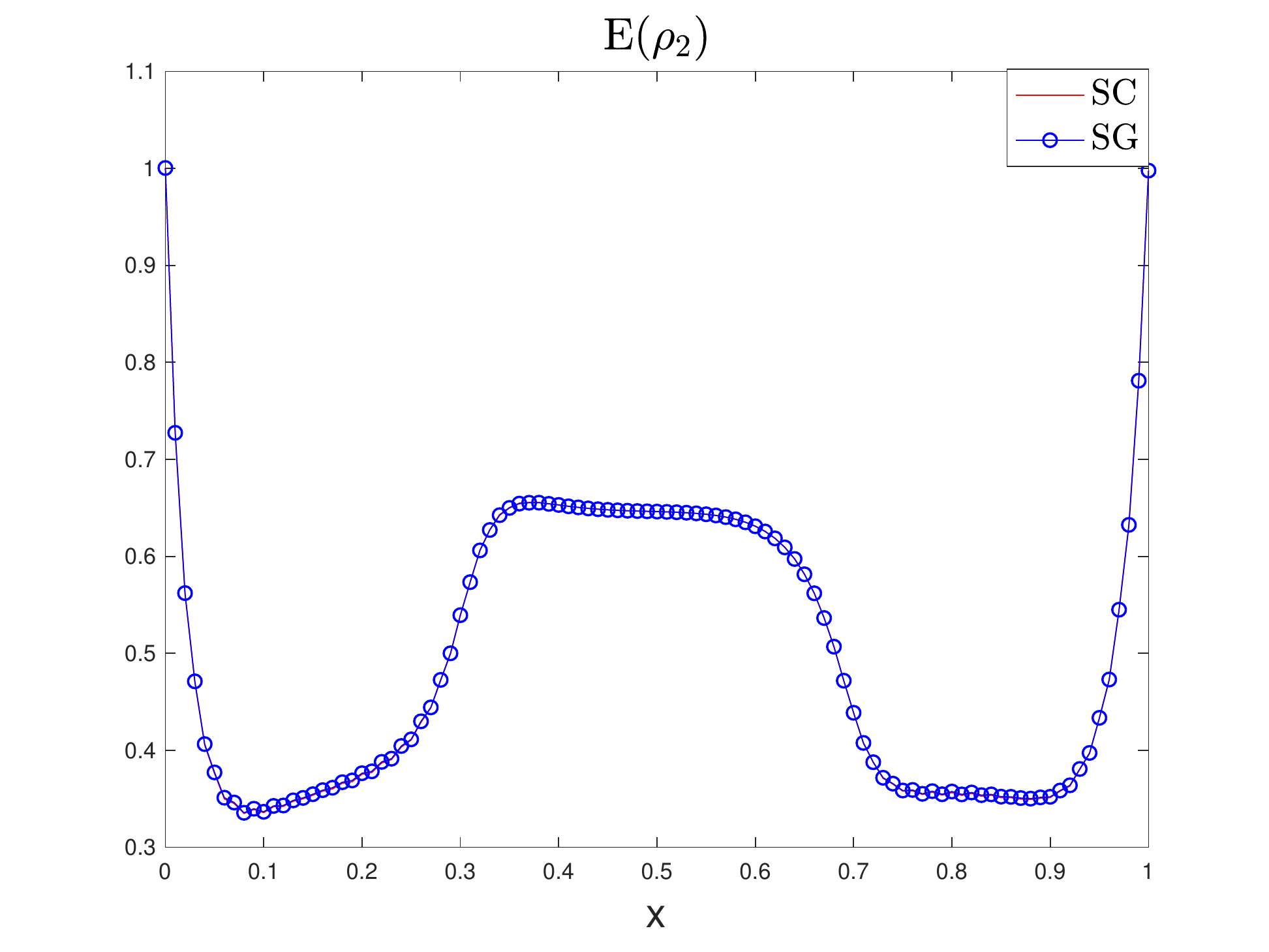}
\centering
\includegraphics[width=0.49\linewidth]{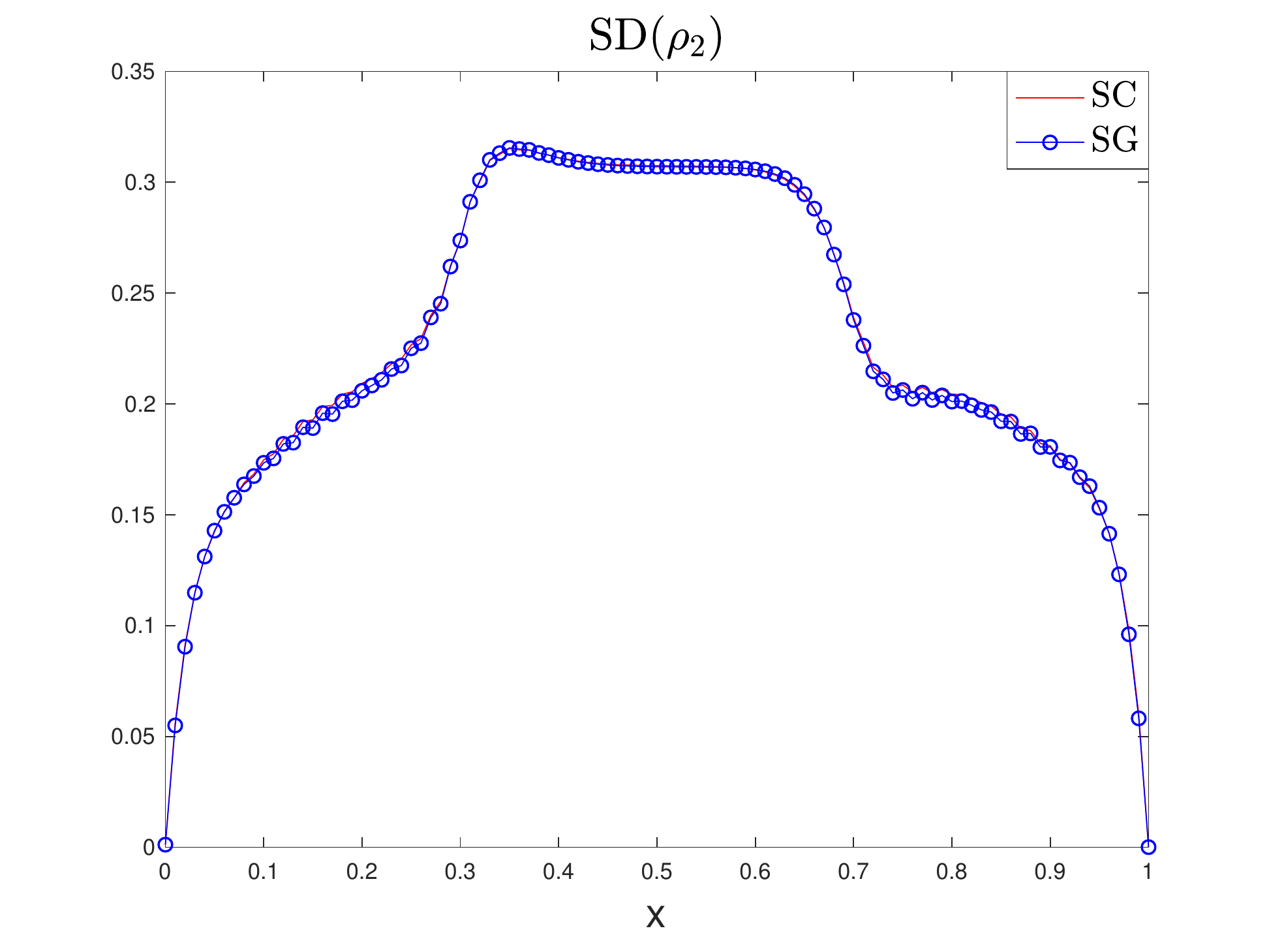}
\caption{Test 2 c).  Red solid line: reference solutions by the SC method. Blue line with circles: gPC-SG method with $K=4$. }
\label{Test2c}
\end{figure}

In Figures \ref{Test2a}, \ref{Test2b} and \ref{Test2c}, 
the high-order stochastic collocation method with $16$ Legendre-Gauss quadrature points is used to obtain the reference solutions. 
A satisfactory agreement between gPC-SG solutions and the reference solutions is clearly observed. 
\\[20pt]
{\bf\large{Test 2 d): Spectral convergence test}}

In this test, the same data and parameters as Test 2 a) are used, where both the doping profile and collision kernels are random. 

\begin{figure} [H]
\centering
\includegraphics[width=0.496\linewidth]{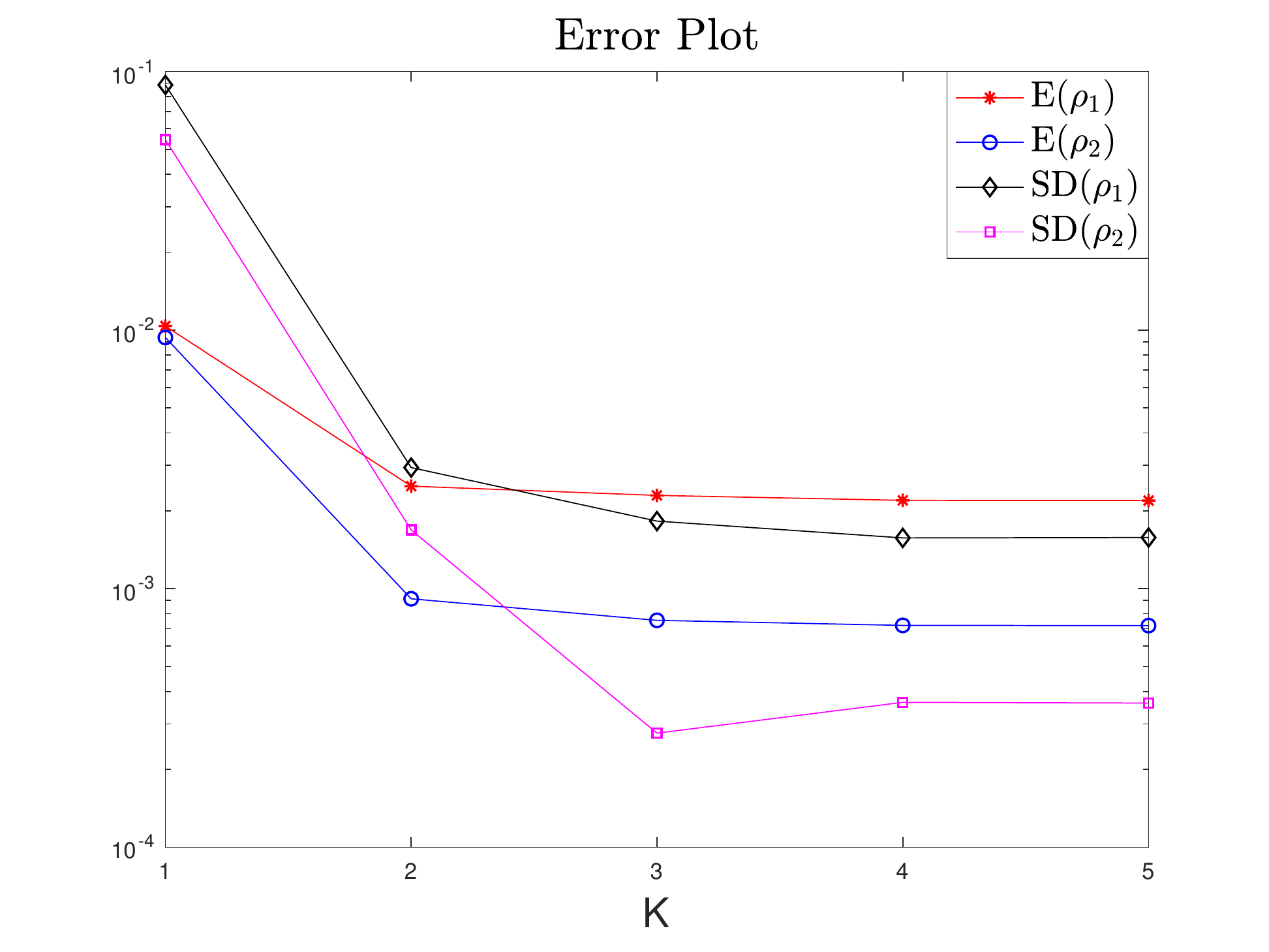}
\centering
 \includegraphics[width=0.496\linewidth]{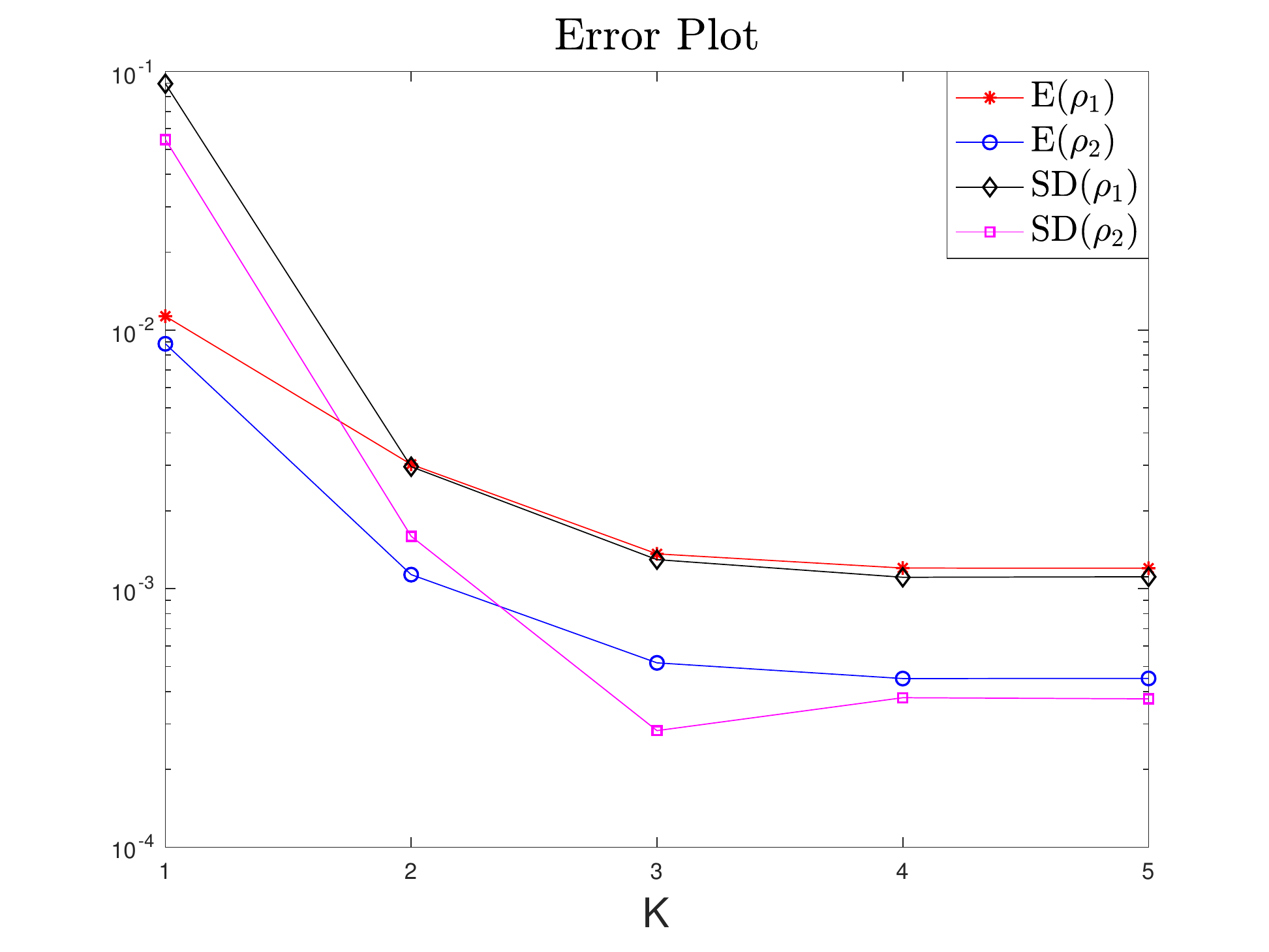}
 \caption{Test 2 d).  Error plots for mean and standard deviation of $\rho_1$, $\rho_2$, 
 $\varepsilon=10^{-3}$ (left) and $\varepsilon=10^{-4}$ (right). Output time is $T=0.005$.  }
 \label{Error_K}
 \end{figure}

Figure \ref{Error_K} shows a semi-log plot for the errors of mean and variance of physical quantities $\rho_1$, $\rho_2$
(density of electrons and holes) 
with $\varepsilon=10^{-3}$ or $\varepsilon=10^{-4}$, using different gPC orders $K$. Error plot for mean and variance of the momentum
give similar results, and we omit it here. We demonstrate a fast exponential convergence with respect to an increasing $K$. 
The errors quickly saturate at modest gPC order $K=4$, then the errors from the temporal and spatial discretization start to dominate and contribute more than 
that from the gPC expansion. This result verifies the s-AP property indicating one can choose $K$ independent of $\varepsilon$. 


\section{Conclusions}
\label{sec:7}

In this paper, we study the bipolar Boltzmann-Poisson model, both for the deterministic problem and the problem with uncertainties, 
with asymptotic behavior leading to the drift diffusion-Poisson system as the Knudsen number goes to zero. 
A s-AP scheme in the gPC-SG framework for the bipolar model with random inputs is designed and numerically verified its efficiency and accuracy. 
Using hypocoercivity of kinetic operators, 
we conduct a convergence rate analysis for both the analytical solution and the gPC solution for a simpler model (without electric field), and conclude their convergence rate exponentially decaying in time, under suitable assumptions. 
A formal proof of s-AP property and a {\it uniform} spectral convergence in the random space for the gPC-SG scheme is obtained. 

{\color{red}Overall, the author thinks that the development of stochastic asymptotic-preserving methods for the bipolar system with random inputs, 
combined with the sensitivity analysis and uniform spectral convergence with an exponential decay in time of the numerical error of the gPC-SG scheme in this project is a {\it first, new and nontrivial} contribution to this field of interest and important for potential applications. }

Future work include conducting a convergence rate analysis for the full model (with the self-consistent electric field); designing and implementing
AP schemes that describe the dynamics of a disparate mass binary gas or plasma system, at various time scales, based on the analysis conducted by Degond and Lucquin-Desreux in \cite{Degond1, Degond2}. 
{\color{blue} Here, we use a second order space discretization and a first order time splitting, similar to that proposed in \cite{JinLorenzo, JPT2}. 
It would be nice to improve the first order time approximation and develop a fully second order scheme, for example, by adopting the method 
introduced in \cite{GL}. This is also considered as a future work. }

\section*{Acknowledgement}

The author would like to thank Prof. Shi Jin and Prof. Irene Gamba for bringing the author's attention 
to work on this project and discussions. {\color{magenta} The author also appreciates both referees' comments on helping improve 
the quality of this paper. }


\newpage 
\section*{Appendix}
\renewcommand{\theequation}{A.\arabic{equation}}

(i) We first show the following two equations needed when deriving the system (\ref{EO1}) from (\ref{EO}): 
\begin{equation}\label{arg1}2\int \sigma(v,w)r(w)dw=\int \sigma(v,w)f(w)dw+\int \sigma(-v,w)f(w)dw, 
\end{equation}
and \begin{equation}\label{arg2}\int\sigma(v,w)j(w)dw=\frac{1}{2\epsilon}\left[\int \sigma(v,w)f(w)dw-\int \sigma(-v,w)f(w)dw\right], 
\end{equation}
for $v>0$. 

Denote $R(v)=\int \sigma(v,w)r(w)dw$, then 
\begin{align}
&\displaystyle R(v)=\int_{w>0}\sigma(v,w)r(w)dw+\int_{w<0}\sigma(v,w)r(w)dw 
=\int_{w>0}\sigma(v,w)r(w)dw+\int_{w>0}\sigma(v,-w)r(w)dw  \notag\\[4pt]
&\label{R}\displaystyle\qquad=\frac{1}{2}\int_{w>0}\sigma(v,w)\left[f(w)+f(-w)\right]dw+\frac{1}{2}\int_{w>0}\sigma(v,-w)\left[f(w)+f(-w)\right]dw. 
\end{align}
For $v>0$, RHS of (\ref{arg1}) is given by
\begin{align}
&\displaystyle\qquad\int\sigma(v,w)f(w)dw+\int\sigma(-v,w)f(w)dw \notag\\[4pt]
&\displaystyle\label{RHS1}=\int_{w>0}\sigma(v,w)f(w)dw+\int_{w<0}\sigma(v,w)f(w)dw
+\int_{w>0}\sigma(-v,w)f(w)dw+\int_{w<0}\sigma(-v,w)f(w)dw  \\[4pt]
&\displaystyle\label{RHS2}=\int_{w>0}\sigma(v,w)f(w)dw+\int_{w>0}\sigma(v,-w)f(-w)dw 
+\int_{w>0}\sigma(v,-w)f(w)dw+\int_{w>0}\sigma(v,w)f(-w)dw \\[4pt]
&\displaystyle \notag = 2 R(v). 
\end{align}
The last step is obvious from (\ref{R}). To check the second equality, we use the 
change of variable $w'=-w$; rotationally invariance and the symmetry of $\sigma$. 

\noindent The third term of (\ref{RHS1}) equals to 
\begin{align*}
&\qquad\int_{w>0}\sigma(-v,w)f(w)dw=\int_{w'<0}\sigma(-v,-w')f(-w')dw'=\int_{w'<0}\sigma(v,w')f(-w')dw' \\[2pt]
&=\int_{w'<0}\sigma(w',v)f(-w')dw'=\int_{w>0}\sigma(-w,v)f(w)dw=\int_{w>0}\sigma(v,-w)f(w)dw, 
\end{align*}
which is the third term of (\ref{RHS2}). The fourth term of (\ref{RHS1}) equals to 
$$\int_{w<0}\sigma(-v,w)f(w)dw=\int_{w'>0}\sigma(-v,-w')f(-w')dw'=\int_{w>0}\sigma(v,w)f(-w)dw, $$
which is the fourth term of (\ref{RHS2}). It is obvious that the first and second term of (\ref{RHS1}) equal to (\ref{RHS2}), respectively. 
Thus we proved (\ref{arg1}). Similarly, one can prove (\ref{arg2}), then we have
$$\int\sigma(v,w)j(w)dw=0, $$ due to the odd function $j$.
\\[10pt]

(ii) We now derive the definitions for the operators $I_{\text{i,plus}}$, $ I_{\text{i,minus}}$. For $v>0$, one has
\begin{align}
&\displaystyle \label{eqn0}\quad\frac{1}{2} \left[I_1(f_1,f_2)(v)+I_1(f_1,f_2)(-v)\right] \\[4pt]
&\displaystyle \label{eqn1}=\frac{1}{2}\int \left( \sigma_I(v,w) + \sigma_I(-v,w)\right)dw\, M_1(v) \\[4pt]
&\displaystyle\quad-\int \sigma_I(v,w)r_2(w)M_2(w)dw\, r_1(v) -\epsilon \int \sigma_I(v,w)j_2(w)M_2(w)dw\, j_1(v)  \notag\\[4pt]
&\displaystyle =\frac{1}{2}\int \left( \sigma_I(v,w) + \sigma_I(-v,w)\right)dw\, M_1(v)-\int \sigma_I(v,w)r_2(w)M_2(w)dw\, r_1(v) \notag\\[4pt]
&\displaystyle \notag := I_{\text{1,plus}}(r_1,r_2),
\end{align}
where $j$ being an odd function is used in the second equality. 
To derive (\ref{eqn1}) from (\ref{eqn0}), note that 
\begin{align*}
&\displaystyle\int \sigma_I(v,w)f_2(w)M_2(w)dw\, f_1(v) + \int\sigma_I(-v,w)f_2(w)M_2(w)dw\, f_1(-v) \\[4pt]
&\displaystyle=\left(\int_{w>0}\sigma_I(v,w)f_2(w)M_2(w)dw +\int_{w>0}\sigma_I(v,-w)f_2(-w)M_2(w)dw \right)
\left(r_1(v)+\epsilon j_1(v)\right) \\[4pt]
&\displaystyle\quad+ \left(\int_{w>0}\sigma_I(v,-w)f_2(w)M_2(w)dw +\int_{w>0}\sigma_I(v,w)f_2(-w)M_2(w)dw\right)
\left(r_1(v)-\epsilon j_1(v)\right), 
\end{align*}
and 
\begin{align*}
&\quad\displaystyle\int\sigma_I(v,w)r_2(w)M_2(w)dw \\[4pt]
&=\displaystyle\frac{1}{2}\int_{w>0}\sigma_I(v,w)\left(f_2(w)+f_2(-w)\right)M_2(w)dw
+\frac{1}{2}\int_{w>0}\sigma_I(v,-w)\left(f_2(w)+f_2(-w)\right)M_2(w)dw, 
\end{align*}
and also 
\begin{align*}
&\quad\displaystyle\int\sigma_I(v,w)j_2(w)M_2(w)dw \\[4pt]
&=\displaystyle\frac{1}{2}\int_{w>0}\sigma_I(v,w)\left(f_2(w)-f_2(-w)\right)M_2(w)dw
-\frac{1}{2}\int_{w>0}\sigma_I(v,-w)\left(f_2(w)-f_2(-w)\right)M_2(w)dw,
\end{align*}
thus it is easy to see that (\ref{eqn1}) equals to (\ref{eqn0}). 
We derived the definition for $I_{1,\text{plus}}$, which can be written as a function of $r_1$ and $r_2$. 

Similarly for $I_{\text{1,minus}}$, one gets
\begin{align}
&\displaystyle\quad\frac{1}{2} \left[I_1(f_1,f_2)(v)-I_1(f_1,f_2)(-v)\right] \notag \\[4pt]
&\displaystyle=\frac{1}{2}\int \left(\sigma_I(v,w) -\sigma_I(-v,w)\right)dw M_1(v) \notag\\[4pt]
&\displaystyle\quad-\epsilon\int \sigma_I(v,w)j_2(w)M_2(w)dw r_1(v) - \epsilon \int \sigma_I(v,w)r_2(w)M_2(w)dw\, j_1(v) \notag \\[4pt]
&\displaystyle=\frac{1}{2}\int \left(\sigma_I(v,w) -\sigma_I(-v,w)\right)dw M_1(v)- \epsilon \int \sigma_I(v,w)r_2(w)M_2(w)dw\, j_1(v) \notag\\[4pt]
&\displaystyle\notag:= I_{\text{1,minus}}(r_2, j_1). 
\end{align}
$I_{\text{2,plus}}$, $I_{\text{2,minus}}$ can be similarly obtained, and we omit the details. 
The definitions of these four operators are given in equations (\ref{I_DEF}). 

\bibliographystyle{siam}
\bibliography{Two_band.bib} 
\end{document}